\documentclass [11 pt , a4 paper , oneside]{article}

\usepackage{color}
\usepackage{amsthm}
\usepackage{subfigure}
\usepackage{amsmath}
\usepackage{mathtools}
\usepackage{authblk}
\usepackage[top=1in, bottom=1.25in, left=1.25in, right=1.25in]{geometry}

\usepackage{tikz}
\usetikzlibrary{chains,positioning,shapes.symbols}
\usetikzlibrary{shapes.geometric, arrows}
\usetikzlibrary{shapes,arrows}
\usetikzlibrary{shapes.callouts}
\usetikzlibrary{positioning}
\usepackage{tikzsymbols}

\input{Pack.sty}

\definecolor{Maria}{HTML}{008b8b} %039141}

\newcommand{\fb}[1]{{\color{black}#1}\color{black}{}}
\newcommand{\B}[1]{{\color{black}#1}\color{black}{}}
\newcommand{\A}[1]{{\color{black}#1}\color{black}{}}

\begin{document}

\title{POD-Galerkin Model Order Reduction for Parametrized Nonlinear Time Dependent Optimal Flow Control: an Application to Shallow Water Equations.}
\author[$\sharp$]{Maria Strazzullo}
\author[$\sharp$]{Francesco Ballarin}
\author[$\sharp$]{Gianluigi Rozza}
\affil[$\sharp$]{mathlab, Mathematics Area, International School for Advanced Studies (SISSA), Via Bonomea 265, I-34136 Trieste, Italy. 

\textnormal{maria.strazzullo@sissa.it, francesco.ballarin@sissa.it, gianluigi.rozza@sissa.it}}
%\affil[$\sharp$]{}
\date{}                     %% if you don't need date to appear
\setcounter{Maxaffil}{0}
\renewcommand\Affilfont{\itshape\small}
\maketitle

\tikzstyle{startstop} = [rectangle, rounded corners, minimum width=3cm, minimum height=1.2cm, text centered, text width=3cm, draw=black, fill=blue!15]
\tikzstyle{io} = [trapezium, trapezium left angle=70, trapezium right angle=110, minimum width=.8cm, minimum height=1.1cm, text centered,  text width=1.5cm, draw=black, fill=red!40]
\tikzstyle{arrow} = [thick,->,>=stealth]
\tikzstyle{b_arrow_dashed} = [dashed,<->,>=stealth]
\tikzstyle{b_arrow} = [thick,<->,>=stealth]
\tikzstyle{line} = [thick,-,>=stealth]
\tikzstyle{line_dot} = [dotted,-,>=stealth]

\begin{abstract}
\no In this work we propose reduced order methods as a reliable strategy to efficiently solve parametrized optimal control problems governed by shallow waters equations in a solution tracking setting. The physical parametrized model we deal with is nonlinear and time dependent: this leads to very time consuming simulations which can be unbearable \fb{e.g.} in a marine environmental monitoring plan application. Our aim is to show how reduced order modelling could help in \fb{studying} different configurations and phenomena in a fast way. After building the optimality system, we rely on a POD-Galerkin reduction in order to solve the optimal control problem in a low dimensional reduced space. The presented theoretical framework is actually suited to general nonlinear time dependent optimal control problems. The proposed methodology is finally tested with a numerical experiment: the reduced optimal control problem governed by shallow waters equations reproduces the desired velocity and height profiles \fb{faster than} the standard model, still remaining accurate.
\end{abstract}

\section{Introduction}
\label{intro}
Parametrized optimal control problems (\ocp s) governed by parametrized partial differential equations (PDE($\boldsymbol{\mu}$)s) are very powerful mathematical formulations, to be exploited in several applications in different scientific fields, see \cite{leugering2014trends} for an overview. Among the possible impacts that \ocp s  can have in scientific research, we will refer to the investigation into problems dealing with environmental sciences. Indeed, this work is motivated by the ongoing demand for reaching fast and accurate simulations for the coastal marine environment safeguard. The marine ecosystem is related to other important social factors such as, for example, economic growth, natural resources preservation, monitoring plans. 
 Furthermore, the marine environment is very far to be completely understood, since it is related to very complicated natural phenomena and anthropic consequences \cite{cavallini2012quasi, mosetti2005innovative, golfo}.
For sure, the parametric setting is necessary in order to study different configurations: the parameter $\boldsymbol{\mu} \in \Cal P \subset \mathbb R^d$ could represent a variety of physical phenomena. Moreover, in the environmental field, the theory of \ocp s fits well with the need of increasing the models forecast capabilities through a data assimilation approach \cite{ghil1991data,kalnay2003atmospheric,tziperman1989optimal}. A lot of effort is made in order to make the predictions of PDEs-based models the most similar to collected data. 
Data assimilation \ocp s have been already analysed in several works, as \cite{phd, quarteroni2005numerical, quarteroni2007reduced, Strazzullo1, ZakiaMaria}. \fb{Yet, }the main drawback of data assimilated problems is the huge computational complexity which still limits their applicability, most of all if the optimization problem deals with very complicated parametric flow models, as the ones used in marine and coastal engineering. Furthermore, in the \fb{described context}, accurate simulations are required in a small amount of time, in order to better study and analyse them rapidly. This is the reason \fb{that motivates 
the use of Reduced Order Methods (ROMs)} as a suitable approach for fast and accurate surrogate simulations of partial differential equations PDE($\boldsymbol{\mu}$)s \cite{hesthaven2015certified, prud2002reliable, RozzaHuynhPatera2008}. The main feature of ROM techniques is to solve the parametrized problem in a low dimensional framework in order to save computational resources which can be exploited for the analysis of several parametric configurations: ROMs recast a time consuming simulation, the \emph{truth problem}, into a new fast and reliable formulation thanks to a Galerkin projection into reduced spaces, generated by basis functions derived from a proper orthogonal decomposition (POD) algorithm, as presented in \cite{ballarin2015supremizer, burkardt2006pod, Chapelle2013, hesthaven2015certified}.
In general, reduction methods for parametrized nonlinear time dependent \ocp s are very complex to analyse both theoretically and numerically. Although the literature is quite consolidated for steady constraints, see for example \cite{bader2016certified,bader2015certified,dede2010reduced, gerner2012certified,karcher2014certified,karcher2018certified,kunisch2008proper,negri2013reduced,negri2015reduced,quarteroni2007reduced},  \fb{where} the interested reader may find theoretical and numerical analysis for different linear models, there is very small knowledge about time dependency \cite{Iapichino2, karcher_grepl_2014, Strazzullo1,ZakiaMaria}. Another difficulty to be overcome is the treatment and the reduction of nonlinear \ocp s, see for example 
\cite{LassilaManzoniQuarteroniRozza2013a, rozza2012reduction, Strazzullo1,ZakiaMaria,Zakia}. 
In this work, we focus on ROM for \ocp s with quadratic cost functional constrained to parametrized Shallow Waters Equations (SWEs). The latter is a very useful model in environmental sciences, which is capable to simulate various marine phenomena such as, for example, tidal flows and mixing, currents action on shorelines and coasts, planetary flows and even tsunamis \cite{cavallini2012quasi,Shallow}. The state equation proposed is nonlinear and time dependent: this leads to growing complexity in the \fb{solution} of the optimality system in a real-time context. Indeed, even if the state equation have been analysed and managed numerically with many approaches, see for example \cite{Agoshkov1993, agoshkov1994recent2, agoshkov1994recent, ferrari2004new, miglio1999finite, miglio2005model, miglio2005model2, ricchiuto2007application, ricchiuto2009stabilized,  takase2010space}, optimal control strategies, see e.g. \cite{agoshkov2007optimal}, and their parametric formulation are still quite unexplored. Moreover reduced techniques has been merely applied to the state equation \cite{Navon2, Navon1,TaddeiSWE}.
\\ \A{The main novelty of this work is to perform reduction on the parametric space on the complete SWEs model, i.e.\ to a nonlinear and time dependent problem, in a solution tracking optimal control framework. Indeed, to the best of our knowledge, there are no contributions on the topic of physical parametrized \ocp s governed by such a model that, despite its complexity, is of growing interest in many fields of applications.}
Thus, we aim at making a further step towards forecasting data assimilated coastal models which could be used as resources to manage realistic experiments involved in marine sciences with environmental prevision purposes.  
\A{In this setting, we want to provide a versatile tool to be exploited in an interdisciplinary framework. As far as we know, reduced algorithm have never been applied to nonlinear time dependent \ocp s. We here propose them to produce a reliable reduced order model that can be effective in providing a large number of parametric simulations in an acceptable amount of time, and this can improve the study of the \fb{considered physical phenomenon}}. \A{
Moreover, we adapt \emph{space-time} reduced techniques, already used for parabolic problems, to nonlinear time dependent \ocp s.
We will show how ROMs could be a good strategy which will give us faster, but still accurate, results, complying with standard techniques already exploited in simpler contexts, such as the steady and linear time dependent governing equations.}
\no The work is outlined as follows. In Section 2, we first present the SWEs model and then we show it in an optimal control framework. Moreover, we briefly describe the discrete approximation and the algebraic version of the presented solution tracking problem. Section 3 will introduce the basic ideas behind ROM discretization for \ocp s \cite{hesthaven2015certified, ito1998reduced, karcher2014certified}: we will describe POD sampling algorithm for \ocp s and the aggregated reduced spaces technique used in \cite{dede2010reduced, negri2015reduced,negri2013reduced} that will guarantee the solvability of the optimality system in its saddle point formulation. Moreover, we will briefly mention the affine decomposition assumption, see e.g.\ \cite{hesthaven2015certified}, needed for an efficient \fb{solve} of the reduced system. In Section 4, we test our methodology on a parametrized \ocp s governed by the SWE equations, inspired by an uncontrolled numerical test case of \cite{ferrari2004new}, where the evolution a Gaussian water height is studied. Our test case aims at recovering a given desired velocity-height profile. Finally, conclusions and perspectives follow in Section 5.

\section{Problem Formulation and Discretization}
In this section, we will introduce a OCP($\bmu$) governed by SWEs and its truth discretization. As already mentioned in Section \ref{intro}, the SWEs are a great tool in order to simulate coastal behaviour. A brief introduction to the state equation follows. Then, the SWEs will be connected to their OCP($\bmu$) framework in Section \ref{ocp_SWEs}. Then, we will describe the full order approximation of our problem based on a Finite Element space-time approach.
\subsection{The Shallow Waters Equations}
\label{SWE}
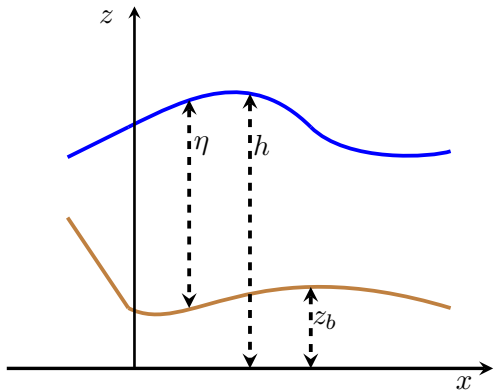
\begin{figure}[H]
\begin{center}
\begin{tikzpicture}[scale=0.8]
%\filldraw [gray] (0,0) circle (2pt) circle (3pt) (2,1) circle (2pt) (2,0) circle (2pt);
%\draw (0,0) .. controls (-0.7,-0.3) and (-1,-.1) .. (-1.5,-.2);
\draw[line width=0.5mm,  blue] (0,0) .. controls (1,0.5) and (2,1) .. (3,0);
\draw [line width=0.5mm,  blue](3,0) .. controls (3.6,-0.6) and (5,-0.5) .. (5.3,-0.4);
\draw[line width=0.5mm,  blue] (0,0) .. controls (-1,-0.5) and (-1,-0.5) .. (-1,-0.5);

\draw[line width=0.5mm,  brown] (0,-3) .. controls (1,-3.5) and (2,-2) .. (5.3,-3);
%\draw [line width=0.5mm,  brown](3,-3) .. controls (3.6,-3.6) and (5,-3.5) .. (5.3,-3.4);
\draw[line width=0.5mm,  brown] (0,-3) .. controls (-1,-1.5) and (-1,-1.5) .. (-1,-1.5);

\draw [arrow, line width=0.35mm] (0.1,-4) -- (0.1, 2);
\draw [arrow, line width=0.45mm] (-2,-4) -- (6, -4);
\draw[b_arrow_dashed,  line width=0.5mm] (1,0.45) -- (1, -3);
\draw[b_arrow_dashed,  line width=0.5mm] (2,0.55) -- (2, -4);
\draw[b_arrow_dashed,  line width=0.5mm] (3,-2.65) -- (3, -4 );
% \draw[line, dashed, line width=0.5mm, gray] (-1,-1) -- (5.3, -1);
\node[anchor=south, label=left:{$\vspace{2mm}x$}] (6,-4.4) at (6,-4.4) {};
\node[anchor=north, label=left:{$\vspace{2mm}z$}] (0.1,2) at (0.1,2) {};
\node[anchor=west, label=left:{$\vspace{2mm}\eta$}] (1.5, -0.3) at (1.5, -0.3) {};
\node[anchor=west, label=left:{$\vspace{2mm}h$}] (2.5,- 0.3) at (2.5,- 0.3) {};
\node[anchor=west, label=left:{$\vspace{2mm}z_b$}]  (3.6,-3.2) at (3.6,-3.2) {};
\end{tikzpicture}
\end{center}
\caption{Notations: schematic representation.}
\label{notations}
\end{figure}

Now we aim at describing the parametrized SWEs. The interested reader may refer to classical references \cite{cavallini2012quasi, Shallow}, where the topic is deeply analysed in its total generality for a space-time domain $Q = \Omega \times (0,T) \subset \mathbb R^2 \times \mathbb R$. This state equation has been studied both from the analytical and numerical point of view in many works, see, for example, \cite{agoshkov1994recent2,  Agoshkov1993, agoshkov1994recent, ferrari2004new, miglio1999finite, miglio2005model, miglio2005model2, ricchiuto2007application, ricchiuto2009stabilized}.\\
Let us define $Y_v = H^1_{\Gamma _{D_v}}(\Omega)$, $Y_h = L^2_{\Gamma_{D_h}}(\Omega)$ and the space  $U = [L^2(\Omega)]^2$, where $\Gamma_{D_v}$ and $\Gamma_{D_h}$ are portions of the boundary domain $\partial \Omega$ where Dirichlet boundary conditions have been imposed for the vertically averaged velocity profile of the wave $\boldsymbol v$ and the free surface elevation variable $h$, respectively. With the term $\eta = h - z_b$ we indicate the water depth, where $z_b$ represents the bottom bathymetry of the domain that we are considering: a schematic description of the involved variables is given in Figure \ref{notations}.  We used the standard 2D-model obtained by the vertical integration of the velocity variable as presented in \cite{Agoshkov1993, miglio1999finite}.
\no
Since we are dealing with time dependent problems, we seek our state solution in the spaces 
$$\Cal Y_v = \Big \{ \boldsymbol v \in L^2(0,T; [Y_v]^2) \text{ such that } \dt{\boldsymbol v} \in L^2(0,T; [Y_v\dual]^2) \Big \},$$ and similarly in
$$\Cal Y_h = \Big \{ h \in L^2(0,T; Y_h) \text{ such that } \dt{h} \in L^2(0,T; Y_h\dual) \Big \},$$
for velocity and height variables, respectively.\\
\no  We will denote \fb{by} $\Cal Y$ the product of the \fb{function} spaces $\Cal Y_v$ and $\Cal Y_h$, i.e. $\Cal Y = \Cal Y_v \times \Cal Y_h$. Moreover, $\Cal Y$ is an Hilbert space with the following norm:
\begin{align*}
\norm{(\star, \cdot)}_{\Cal Y}^2 & =
\norm{\star}_{\Cal Y_{\boldsymbol v}}^2
+ \norm{\cdot}_{\Cal Y_{h}}^2 		\\
& = \norm{\star}_{L^2(0,T; [Y_{\boldsymbol v}]^2)}^2 + \parnorm{\dt{\star}}_{L^2(0,T; [Y_{\boldsymbol v}\dual]^2)}^2
+ \norm{\cdot}_{L^2(0,T; Y_{h})}^2 + \parnorm{\dt{\cdot}}_{L^2(0,T; Y_{h}\dual)}^2.
\end{align*}
\noindent We considered a SWEs model with the following parametrized formulation: given $\boldsymbol \mu \in \Theta \subset \mathbb R^2$ and a forcing term $\boldsymbol u \in \Cal U := L^2(0,T;U)$, find  the parameter dependent \fb{pair} $(\boldsymbol v, h) \in \mathcal Y$ which \fb{satisfies}:
\begin{equation}
\label{SWEs_state}
\begin{cases}
\displaystyle \dt{\boldsymbol v} + \mu_1 \Delta \boldsymbol v + \mu_2(\boldsymbol v \cdot \nabla)\boldsymbol v + g \nabla \eta - \boldsymbol u = 0 & \text{in } Q , \vspace{1mm}\\
\displaystyle \dt{h} + \dive(\eta \boldsymbol v) = 0 & \text{in } Q, \\
\boldsymbol v = \boldsymbol v_0 & \text{on } \Omega \times \{  0 \}, \\
h = h_0 & \text{on }  \Omega \times \{  0 \}, \\
\boldsymbol v = \boldsymbol 0& \text{on }  \partial \Omega \times [0, T].
\end{cases}
\end{equation}
We stress that \fb{the solution depends on the parameter $\bmu = (\mu_1, \mu_2)$, i.e. $(\boldsymbol v, h):= (\boldsymbol v(\bmu), h(\bmu) )$, but in the following we will omit such dependence for compactness of notation}. 
The \fb{proposed analysis} does not change with other general boundary conditions: for an insight on the admissible boundary conditions we refer to \cite{agoshkov1994recent}.
The \fb{proposed model} describes free surface incompressible flows under the assumption of hydrostatic pressure. This hypothesis is valid when the water height is much lower than the wavelength: this is the case of coastal behaviour and shallow depths. In equation \eqref{SWEs_state}, we simplistically represent with the forcing term $\boldsymbol u$ all the physical quantities which can affect the dynamic of the solution, such as the wind stress and the bottom friction: the reason of this choice will be clarified in the next Section. Moreover, we deal with the following parametric context: $\mu_1$ and $\mu_2$ represent how diffusion and advection phenomena affect the shallow waters system, respectively. Furthermore, we assume the bottom to be constant with respect to time and spatial variables for simplicity. Indeed, under this assumption, the bottom has no influence on the system considered, since:
$$
\nabla{\eta} = \nabla{h}  \hspace{2mm} \text{and} \hspace{2mm} \dive{(\eta \boldsymbol v)} = \dive{(h \boldsymbol v)} - z_b\underbrace{\dive{(\boldsymbol v)}}_{= 0} = \dive{(h \boldsymbol v)}.\\
$$
In the next Section, we propose an optimal control problem governed by SWEs, briefly introducing Lagrangian formulation and derivation of the optimality system.
\subsection{Optimal control problem governed by Shallow Waters Equations}
\label{ocp_SWEs}
In this Section we are going to introduce a time dependent  OCP($\bmu$) governed by parametrized SWEs. We will combine the state equation \eqref{SWEs_state} to a minimization problem. The goal is to relate the optimal control theory coupled with SWEs coastal water model as a powerful tool in order to manage marine monitoring issues. Theoretically, we will follow the general theory for time dependent \ocp s proposed in \cite[chapter 3]{troltzsch2010optimal}.
In Section \ref{SWE} we already introduced the \emph{state variable} 
$(\boldsymbol v, h) \in \mathcal Y$. In order to set up an OCP($\bmu$) we need to define a \emph{control variable} $\boldsymbol u$ in the \fb{function} space $\Cal U$. In our applications, we deal with a distributed optimal control, i.e., the control variable represents the forcing term comprising the wind action, atmospheric pressure and the bottom slope effect.  Then, even if we are not actually controlling the system, the optimal control framework can be interpreted as an inverse problem which specifies the physical conditions guaranteeing a desired velocity-height profile $(\boldsymbol v_d, h_d) \in \mathcal Y$. Let us define the \emph{state-control variable} $x = ((\boldsymbol v, h), \boldsymbol u) \in \mathcal X = \mathcal Y \times \mathcal U$. As already stressed in the previous section, even if the variables are $\boldsymbol \mu -$dependent, for the sake of notation we will use 
$((\boldsymbol v, h), \boldsymbol u) := ((\boldsymbol v(\boldsymbol \mu), h(\boldsymbol \mu)), \boldsymbol u(\boldsymbol \mu))$ and 
$x: = x(\boldsymbol \mu)$.
\\
The parametrized OCP($\bmu$) has the following formulation: given a parameter $\boldsymbol \mu \in \Cal P$, find 
$x \in \mathcal X$ which minimized the function $J(x)$ defined as 
\begin{align}
\label{J}
\half \intTimeSpace{(h - h_d(\mu_3))^2} + \half \intTimeSpace{(\boldsymbol v - \boldsymbol v_d(\mu_3))^2} +
\alf \intTimeSpace{\boldsymbol u^2},
\end{align}
under the constraint \eqref{SWEs_state}, where $\alpha \in (0,1]$ is a penalization parameter over the control action. \B{Namely, the optimal control problem depends on a physical parameter $\boldsymbol \mu = (\mu_1, \mu_2, \mu_3)$ in the parameter space $\Cal P \subset \mathbb R^3$. The role of the parameters $\mu_1$ and $\mu_2$ was already introduced in Section \ref{SWE}, while the component $\mu_3$ affects the desired solution profile as one can observe from the functional \eqref{J}. 
}In order to solve the problem, we exploit a Lagrangian approach which allowed us to build the optimality conditions that will be discretized in the following sections: in other words we applied an \emph{optimize-then-discretize} strategy, which first derives the optimality system at the continuous level and approximates it only at the end of the procedure as presented in \cite{fernandez2003control}.  \A{First of all, the problem can be read as:
given a parameter $\boldsymbol \mu \in \Cal P$, find 
$x\in \mathcal X$ which minimizes \eqref{J} such that the weak state equation 
$\Cal S(x, (\boldsymbol \kappa, \xi); \boldsymbol \mu)=0$ is verified for all $(\boldsymbol \kappa, \xi) \in \Cal Y$, where $\Cal S(x, (\boldsymbol \kappa, \xi); \boldsymbol \mu)=0$ will denote the weak formulation of \eqref{SWEs_state}. After defining an \emph{adjoint variable} $(\boldsymbol \chi, \lambda) \in \mathcal Y$, we build the following Lagrangian functional 
\begin{equation}
\label{Lagrangian}
\Lg((\boldsymbol v, h), \boldsymbol u,  (\boldsymbol \chi, \lambda)) = J(x) + \Cal S(x, (\boldsymbol \chi, \lambda); \boldsymbol \mu).
\end{equation}
To perform a constrained minimization of \eqref{J}, we differentiate with respect to the variables $((\boldsymbol v, h), \boldsymbol u, (\boldsymbol \chi, \lambda))$ obtaining the optimality system:
\begin{equation}
\label{optimality_system_SWEs}
\begin{cases}
D_{\boldsymbol v}\Lg((\boldsymbol v, h), \boldsymbol u, (\boldsymbol \chi, \lambda))[\boldsymbol z] = 0 & \forall \boldsymbol z \in \Cal Y_v,\\
D_h\Lg((\boldsymbol v, h), \boldsymbol u, (\boldsymbol \chi, \lambda))[q] = 0 & \forall q \in \Cal Y_h,\\
D_{\boldsymbol u}\Lg((\boldsymbol v, h), \boldsymbol u, (\boldsymbol \chi, \lambda))[\boldsymbol \tau] = 0 & \forall \boldsymbol \tau \in \Cal U,\\
D_{\boldsymbol \chi} \Lg((\boldsymbol v, h), \boldsymbol u, (\boldsymbol \chi, \lambda)) [\boldsymbol \kappa]= 0 & \forall \boldsymbol \kappa \in \Cal Y_v,\\
D_\lambda\Lg((\boldsymbol v, h), \boldsymbol u, (\boldsymbol \chi, \lambda))[\xi] = 0 & \forall \xi \in \Cal Y_h.\\
\end{cases}
\end{equation}
The differentiation with respect to the adjoint variable $(\boldsymbol \chi, \lambda)$ coincides with the \emph{state equation} that in strong form reads as \eqref{SWEs_state}. Moreover, differentiation with respect to the control variable $\boldsymbol u$ leads to the \emph{optimality equation} in $Q$ of the form
$
\alpha \boldsymbol u = \boldsymbol \chi,
$ while we can derive the \emph{adjoint equation} differentiating with respect to the state variable and its strong form is 
\begin{align*}
\begin{cases}
\displaystyle \boldsymbol v - \dt{\boldsymbol \chi} + \mu_1 \Delta \boldsymbol \chi - \mu_2 (\boldsymbol v \cdot \nabla)\boldsymbol \chi + \mu_2 (\nabla \boldsymbol v)^T\boldsymbol \chi  - h\nabla \lambda = \boldsymbol v_d & \text{in } Q, \vspace{1mm}\\
\displaystyle  h - \dt{\lambda} - \boldsymbol v \cdot \nabla \lambda - g\dive(\boldsymbol \chi) = h_d & \text{in } Q, \vspace{1mm}\\
\boldsymbol \chi = \boldsymbol 0  & \text{ on } \partial \Omega \times (0,T) \\
\boldsymbol \chi = \boldsymbol 0 & \text{on } \Omega \times \{  T \}, \\
\lambda= 0 & \text{on }  \Omega \times \{  T \}. \\
\end{cases}
\end{align*}
The three equations combined together will denote the strong form of the optimality system we will deal with and that will be solved through numerical approximation. Some details on the weak formulation of the problem can be found in Appendix \ref{1}. In the next Section we propose a space-time approach as full order numerical discretization in a all-at-once framework, as presented in \cite{HinzeStokes, HinzeNS, Stoll1, Stoll} for linear and nonlinear constraints.}

\subsection{Space-Time  \ocp s: All-at-Once Approach}
\label{FEM}
\A{
In this Section, we will present the space-time discretization of the OCP($\bmu$) defined in Section \ref{ocp_SWEs}.
Indeed, a first discretization is needed for the Reduced Order  Model (ROM) approximation, as we will clarify later in Section \ref{sec_ROM}. Namely,  a \textit{truth} problem \fb{solution} is a necessary step in order to build reduced basis functions to apply model reduction. \\
 First of all, we will focus on the Finite Element (FE) approximation. 
As already introduced, our aim is to build a discretized optimality system based on the \emph{first optimize, then discretize} approach, see e.g. \cite{fernandez2003control}. Namely, we first derive the optimality conditions and then we perform a discretization in time and space through Euler's methods and FE approximation, respectively. \\
For this purpose, we define a triangulation $\Cal T$ over the spatial domain $\Omega$. We can now provide the FE spaces as $Y^{\Cal N_{\boldsymbol v}}_v = [Y_v]^2 \cap\mathscr X_{r_{\boldsymbol v}}$, $Y^{\Cal N_{h}}_h = Y_h \cap\mathscr X_{r_h}$
and $U^{\Cal N_{\boldsymbol u}}= U  \cap \mathscr X_{r_{\boldsymbol u}}$, where
$
 \mathscr X_r = \{ s \disc \in C^0(\overline \Omega) \; : \; s |_{K} \in \mbb P^r, \; \; \forall K \in \Cal T \disc \}. 
$ 
The space $\mbb P^r$ consists of all the polynomials of degree at most equal to $r$ and $K$ is a triangular element of $\Cal T$. Let us refer to $\Cal N$ as the global FE dimension of the system, i.e. $\Cal N = 2 \Cal N_{\boldsymbol v} + 2 \Cal N_{h} + \Cal N_{\boldsymbol u}$. Indeed, in this new configuration, the state and adjoint velocity belong to 
$$
\Cal Y_v\disc = \Big \{ \boldsymbol v \in L^2(0,T; Y^{\Cal N_{\boldsymbol v}}_v) \text{ such that } \dt{\boldsymbol v} \in L^2(0,T; {Y^{\Cal N_{\boldsymbol v}}_v}^{\ast}) \Big \},
$$ and, similarly, the state and adjoint elevation variables are in 
the space 
$$
\Cal Y_h\disc = \Big \{ h \in L^2(0,T; Y^{\Cal N_h}_h) \text{ such that } \dt{h} \in L^2(0,T; {Y^{\Cal N_h}_h}^{\ast}) \Big \} .$$
Finally, the function space considered for state and adjoint velocity-height variables is $\Cal Y \disc=  \Cal Y_v\disc \times \Cal Y_h\disc$, while the control variable is in $\Cal U \disc = L^2(0,T; U^{\Cal N_{\boldsymbol u}})$. For the sake of notation, we used as apex the global dimension $\Cal N$ over the spaces, even if it is not the actual dimension of the space considered. 
Indeed, we are dealing with finite dimensional Hilbert spaces and with the basis functions $(\{ \boldsymbol \phi^i \}_{i =1}^{\Cal N_{\boldsymbol v}},\{ \phi^i \}_{i =1}^{\Cal N_{h}})$ and $\{ \boldsymbol \psi^i \}_{i =1}^{\Cal N_{\boldsymbol u}}$ for state/adjoint and control spaces, respectively. \\
 \no The parametrized FE optimality system reads: given $\boldsymbol \mu \in \Cal P$, find the discrete variable $((\boldsymbol v\disc, h\disc), \boldsymbol u\disc, (\boldsymbol \chi \disc, \lambda \disc))$ which solves:
\begin{equation}
\label{FE_optimality_system_SWEs}
\begin{cases}
D_{\boldsymbol v}\Lg((\boldsymbol v \disc, h\disc), \boldsymbol u \disc, (\boldsymbol \chi \disc, \lambda \disc))[\boldsymbol z] = 0 & \forall \boldsymbol z \in \Cal Y_v\disc,\\
D_h\Lg((\boldsymbol v\disc, h\disc), \boldsymbol u\disc, (\boldsymbol \chi \disc, \lambda \disc))[q] = 0 & \forall q \in \Cal Y_h\disc,\\
D_{\boldsymbol u}\Lg((\boldsymbol v\disc, h\disc), \boldsymbol u\disc, (\boldsymbol \chi \disc, \lambda \disc))[\boldsymbol \tau] = 0 & \forall \boldsymbol \tau \in \Cal U \disc,\\
D_{\boldsymbol \chi} \Lg((\boldsymbol v\disc, h\disc), \boldsymbol u\disc, (\boldsymbol \chi \disc, \lambda \disc)) [\boldsymbol \kappa]= 0 & \forall \boldsymbol \kappa \in \Cal Y_v \disc,\\
D_\lambda\Lg((\boldsymbol v\disc, h\disc), \boldsymbol u\disc, (\boldsymbol \chi \disc, \lambda \disc))[\xi] = 0 & \forall \xi \in \Cal Y_h \disc .\\
\end{cases}
\end{equation}
\no As we did in Section \ref{ocp_SWEs}, we indicate the state-control variable $((\boldsymbol v\disc, h\disc),  \boldsymbol u\disc)$ with $x\disc$. We now deal with the time approximation.
The time interval $(0,T)$ is divided in $N_t$ equispaced subintervals with $\Delta t$ as time step. Indeed, at each time $t_k = k \times \Delta t$ for $k= 1, \cdots, N_t$, our FE solution variables can be respectively written as:
\begin{equation*}
\boldsymbol v_k \disc = \sum_{1}^{\Cal N_{\boldsymbol v}}v^{i}_k\boldsymbol \phi^{i}, \hspace{3mm}
h_k \disc = \sum_{1}^{\Cal N_h}h^{i}_k\phi^{i}, \hspace{3mm}
\boldsymbol u_k \disc= \sum_{1}^{\Cal N_{\boldsymbol u}}u^{i}_k\boldsymbol \psi^{i}, \hspace{3mm}
\end{equation*}
\begin{equation*}
\boldsymbol \chi_k \disc = \sum_{1}^{\Cal N_{\boldsymbol v}}\chi^{i}_k\boldsymbol \phi^{i}, \hspace{3mm} \text{and} \hspace{3mm}
\lambda_k \disc = \sum_{1}^{\Cal N_h}\lambda^{i}_k\phi^{i}.
\end{equation*}
Following the strategy presented in \cite{HinzeStokes, Stoll1, Stoll} for linear state equations and in \cite{HinzeNS} for Navier-Stokes equations, we define
$\bar v = [v_1, \dots, v_{N_t}]^T$, $\bar h = [h_1, \dots, h_{N_t}]^T$ and $\bar u = [u_1, \dots, u_{N_t}]^T$, where  $v_k, h_k$ and $u_k$ are the row vectors of the FE coefficients for state discrete variables at each time step. The vectors representing the initial condition for the velocity field $\boldsymbol v$ and the water height $h$ are 
$\bar v_0 = [v_0, 0, \dots, 0]^T$ and $\bar h_0 = [h_0, 0, \dots, 0]^T$, respectively. Following the same argument, we can define both adjoint vectors $\bar \chi = [\chi_1, \dots, \chi_{N_{t}}]^T$ and $\bar \lambda = [\lambda_1, \dots, \lambda_{N_{t}}]^T$ and the desired profiles $\bar v_d = [{v_d}_1, \dots, {v_d}_{N_{t}}]^T$ and $\overline h_d = [{h_d}_1, \dots, {h_d}_{N_{t}}]^T$. Also in this case, $\chi_k$ and $\lambda_k$ represent the vectors of the component of the FE variables at the $k-$th time step, for $k = 1, \dots, N_{t}$. Moreover, from now on, with 
$\boldsymbol w \disc := (\boldsymbol v \disc, h \disc, \boldsymbol u \disc, \boldsymbol \chi \disc, \lambda \disc)$ we refer to the global FE variable including all the time instances, i.e. the $\bigstar \disc$ for $\bigstar = \boldsymbol v, h, \boldsymbol u, \boldsymbol \chi, \lambda$ will indicate the global discretized space-time variable and $\bigstar_k \disc$ will be the global variable evaluated at the time $t_k$.
 The shown structure is consistent with the space-time formulation exploited in several works as \cite{urban2012new, yano2014space, yano2014space1}: in this specific case, the backward Euler scheme in time coincides  with a piecewise constants Discontinuous Galerkin approach, as underlined in \cite{eriksson1987error}.  Although, for the sake of simplicity, we will always refer to Euler's schemes. \\
First of all, let us proceed with the discretization of the \emph{state equation} governing the problem \eqref{SWEs_state}. Using a backward Euler, the \emph{state equation} is discretized forward in time.  The discretization gives the following result for the governing equation at each time step, for $k \in  \{0, \dots, N_t - 1\}$:
\begin{equation}
\label{disc_state} 
\begin{cases}
\displaystyle \frac{\boldsymbol v_{k+1} \disc - \boldsymbol v_k \disc }{\Delta t} + \mu_1 \Delta \boldsymbol v_{k+1} + \mu_2(\boldsymbol v_{k+1} \disc \cdot \nabla) \boldsymbol v_{k+1} \disc + g \nabla (h_{k+1}\disc) = \boldsymbol u_{k+1} \disc \\
\displaystyle \frac{h_{k+1}\disc - h_{k}\disc}{\Delta t} + \dive(h_{k+1} \disc \boldsymbol v_{k+1}\disc) = 0 \\
\end{cases}
\end{equation}
The same discretization strategy can be applied for the \emph{optimality equation}. In this case, at each time step, one solves the equation:
\begin{align}
\label{opt_eq_k}
\alpha \Delta t \boldsymbol u_k\disc = \Delta t \boldsymbol \chi_k \disc& \spazio \text{for } k \in  \{1, \dots, N_t \}.
\end{align}
The last step of the full order discretization involves the \emph{adjoint equation}. Since we are given the value of the adjoint variables at time $T$, the equations are discretized backward in time, through a forward Euler's method. In this case, at each time step, we have to solve the following system for $k \in  \{N_t, \dots, 2\}$:
\begin{equation}
\label{disc_adj} 
\begin{cases}
\displaystyle \boldsymbol v_{k-1}\disc + \frac{\boldsymbol \chi_{k-1}\disc - \boldsymbol \chi_{k}\disc}{\Delta t} + \mu_1 \Delta \boldsymbol \chi_{k - 1}\disc - \mu_2 (\boldsymbol v_{k - 1}\disc \cdot \nabla)\boldsymbol \chi_{k - 1}\disc 
\\ \qquad \qquad \qquad \qquad + \mu_2 (\nabla \boldsymbol v\disc_{k - 1})^T\boldsymbol \chi_{k - 1}\disc  - h_{k - 1}\disc \nabla \lambda_{k - 1}\disc = {\boldsymbol {v}_d}_{k - 1}\disc,\\
\displaystyle  h_{k - 1}\disc + \frac{ \lambda_{k-1}\disc - \lambda_{k}\disc}{\Delta t} - \boldsymbol v_{k - 1}\disc \cdot \nabla \lambda_{k - 1}\disc - g\dive(\boldsymbol \chi_{k - 1} \disc) = {{{h_d}^{\Cal N}}_{k - 1}}.\\
\end{cases}
\end{equation}
We now have all the ingredients to define the whole discrete optimality system, i.e. given a $\boldsymbol \mu \in \Cal P$ find the vector 
$$\bar w = 
\begin{bmatrix}
\begin{bmatrix}\bar v \\ \bar h \\ 
\end{bmatrix}
\\
 \bar u \\
\begin{bmatrix}\bar \chi \\ \bar \lambda \\ 
\end{bmatrix}
\end{bmatrix},
$$ which solves the following nonlinear system
\begin{equation}
\label{G_compact}
\Cal R (\boldsymbol w \disc, \bmu) = \fb{\Cal G}(\boldsymbol w \disc; \boldsymbol \mu) \bar w - \bar f = 0,
\end{equation}
where $\Cal R(\bar w, \bmu)$ represents the residual vector given by the difference of the action the aforementioned nonlinear optimality equations in matrix form and the right hand side vector, respectively denoted with  $\fb{\Cal G}(\boldsymbol w \disc, \boldsymbol \mu)$ and $\bar f$. 
In order to find the space-time optimal solution $\bar w$, we rely on Newton's method, i.e., defining 
$\mathbb J(\bar w; \bmu) = \fb{\mathbb{D}({\Cal G}(\boldsymbol w\disc; \bmu)}\bar w)$ the Frech\'et  derivative of the operator $\fb{\Cal G}(\boldsymbol w \disc; \boldsymbol \mu)\bar w$, we iterate the solution
\begin{equation}
\bar {w}^{j + 1} := \bar {w}^j+ \mathbb J(\bar w^{j}; \boldsymbol \mu)^{-1}(- \Cal R({\boldsymbol w \disc}^j; \boldsymbol \mu)), \spazio j \in \mathbb N,
\end{equation}
until the convergence is reached. \\
In Appendix \ref{2}, we provide the formulation of all the involved discrete quantities. The linearized optimality system carries out a saddle point structure, i.e.\
\begin{equation}
\label{Frechet}
\mathbb J (\bar w; \boldsymbol \mu) \bar w = 
\begin{bmatrix}
\mathsf A & \mathsf B^T \\
\mathsf B & 0 \\
\end{bmatrix}
\begin{bmatrix}
\bar x \\
\bar p
\end{bmatrix}.
\end{equation}
For some properly defined matrices $\mathsf A$ and $\mathsf B$ and vectors $\bar x$ and $\bar p$ (see Appendix \ref{2}). We here underline this peculiar structure since it will help to understand some important concepts of the next Section.  As already specified in Section \ref{intro}, in a parametric context, space-time solutions could be unfeasible due to the huge computational effort required since the system to be solved has $\Cal N_{\text{tot}} = N_t \times \Cal N$ as actual dimension. In the next Section, we will describe reduced order modelling (ROM) techniques, that we use in order to overcome the problem of finding the parametric \fb{solution} of an expensive optimal control system.
}

\section{ROM Approximation for \ocp s}
\label{sec_ROM} 
In this Section, we provide a brief introduction on ROM approximation techniques and we show how to exploit it in the \fb{solution} of SWEs optimal control parametrized systems. Even if we propose the reduced strategy for a very specific governing state equation, the approach could be used for general problems: indeed, we refer to \cite{Strazzullo1, Zakia} for previous applications to steady nonlinear \ocp s. Besides, in order to deal with time dependency, we follow the numerical strategy already presented in \cite{Strazzullo2, ZakiaMaria}. We start with some basic ideas that guarantee an efficient ROM applicability and then, in Sections \ref{POD} and \ref{aggr}, we will move towards the space-time Proper Orthogonal Decomposition (POD) algorithm, see for example \cite{Strazzullo2,ZakiaMaria} as references, as an adaptation of what is already known for linear \ocp s.
\subsection{Reduced Problem Formulation}
In Section \ref{SWE}, we proposed optimal flow control as a way to formulate inverse problems in marine environmental sciences, exploiting the velocity-height model of SWE. As we already specified in Section \ref{intro}, the space-time method could be computationally unfeasible \fb{when interested in solving} several instances of the proposed \ocp, most of all, in a parametrized setting. ROM techniques replace the \emph{truth} system, with a surrogate one, which is \fb{often} smaller \fb{in terms of} dimension.
We now briefly introduce ROM ideas in the \ocp s setting. In order to clarify the role of the parametric setting, we will explicit the $\boldsymbol \mu -$dependency in the quantities that are involved in the reduction process. \\
\fb{By varying} the value of $\boldsymbol\mu$ in the parameter space $\Cal P$, the parametric solution of \eqref{optimality_system_SWEs} will define a manifold
$$
\mathscr M = \{ (\boldsymbol v(\boldsymbol{\mu}), h(\boldsymbol{\mu}), \boldsymbol u(\boldsymbol{\mu}),
\boldsymbol \chi (\boldsymbol{\mu}), \lambda (\boldsymbol \mu))\;| \; \boldsymbol{\mu} \in \Cal P\},
$$
\fb{which we assume to be smooth}.
If we restrict our attention to the space-time approximation, the ensemble of the \emph{truth} solutions is an approximation of $\mathscr M$:
$$
\mathscr M^{\Cal N_{\text{tot}}} = \{ (\boldsymbol v\disc(\boldsymbol{\mu}), h\disc(\boldsymbol{\mu}), \boldsymbol u\disc(\boldsymbol{\mu}),
\boldsymbol \chi \disc (\boldsymbol{\mu}), \lambda \disc (\boldsymbol \mu))\;| \; \boldsymbol{\mu} \in \Cal P\}.
$$
\no Also in this case, the variables $\bigstar \disc$ are actually considered in the space-time function spaces of dimension $N_t \times \Cal  N_{\bigstar}$ for $\bigstar = \boldsymbol v, h, \boldsymbol u, \boldsymbol \chi, \lambda$, with $\Cal  N_{\boldsymbol v} = \Cal  N_{\boldsymbol \chi}$ and $\Cal  N_{ h} = \Cal  N_{\lambda}$, since the same discretized space is used for state and adjoint variables.
ROM aims at describing the structure of the approximated solution manifold $\mathscr M^{\Cal N_{\text{tot}}}$ through the construction of bases derived from \emph{snapshots}, i.e. properly chosen space-time solutions of the variables involved in the optimization system. In other words, reduced spaces are subspaces of the full order spaces and they are chosen through algorithms that guarantee a proper description of how the optimality system \eqref{optimality_system_SWEs} changes with respect \fb{to} a new value of $\boldsymbol \mu$. After the basis functions building process, a standard Galerkin projection is performed, in order to find a low-dimensional surrogate solution, which is computationally efficient and still accurate with respect to the previous model. 
Let us assume  to have already built the reduced \fb{function} spaces $\Cal Y_N \subset \Cal Y \disc \subset \Cal Y$ and $\Cal U{_N} \subset \Cal U \disc \subset \Cal U$ for state/adjoint variables and control, respectively. The projected reduced \ocp $\:$ reads: given $\boldsymbol \mu \in \Cal P$, find $((\boldsymbol v_N(\bmu), h_N(\bmu)), \boldsymbol u_N (\bmu), (\boldsymbol \chi_N (\bmu), \lambda_N (\bmu)))$ which solves:
\begin{equation}
\label{ROM_optimality_system_SWE}
\begin{cases}
D_{\boldsymbol v}\Lg((\boldsymbol v_N(\boldsymbol \mu), h_N(\boldsymbol \mu)), \boldsymbol u_N(\boldsymbol \mu), (\boldsymbol \chi_N(\boldsymbol \mu), \lambda_N(\boldsymbol \mu)))[\boldsymbol z] = 0 & \forall \boldsymbol z \in \Cal Y_v{_N},\\
D_h\Lg((\boldsymbol v_N(\boldsymbol \mu), h_N(\boldsymbol \mu)), \boldsymbol u_N(\boldsymbol \mu), (\boldsymbol \chi_N(\boldsymbol \mu), \lambda_N(\boldsymbol \mu)))[q] = 0 & \forall q \in \Cal Y_h{_N},\\
D_{\boldsymbol u}\Lg((\boldsymbol v_N(\boldsymbol \mu), h_N(\boldsymbol \mu)), \boldsymbol u_N(\boldsymbol \mu), (\boldsymbol \chi_N(\boldsymbol \mu), \lambda_N(\boldsymbol \mu)))[\boldsymbol \tau] = 0 & \forall \boldsymbol \tau \in \Cal U_N,\\
D_{\boldsymbol \chi} \Lg((\boldsymbol v_N(\boldsymbol \mu), h_N(\boldsymbol \mu)), \boldsymbol u_N(\boldsymbol \mu), (\boldsymbol \chi_N(\boldsymbol \mu), \lambda_N(\boldsymbol \mu))) [\boldsymbol \kappa]= 0 & \forall \boldsymbol \kappa \in \Cal Y_v{_N},\\
D_\lambda\Lg((\boldsymbol v_N(\boldsymbol \mu), h_N(\boldsymbol \mu)), \boldsymbol u_N(\boldsymbol \mu), (\boldsymbol \chi_N(\boldsymbol \mu), \lambda_N(\boldsymbol \mu)))[\xi] = 0 & \forall \xi \in \Cal Y_h{_N}.\\
\end{cases}
\end{equation}
The reduced system \eqref{ROM_optimality_system_SWE} is still nonlinear and it can be solved thanks to a Newton's method, as already specified in Section \ref{FEM}. In the next Sections, we will show an approach that leads to the construction of the reduced spaces and what are the techniques \fb{to be used} in order to deal with the reduced Frech\'et derivative \fb{aiming at preserving} the saddle point structure shown in \eqref{Frechet} and its numerical stability.
\subsection{POD Algorithm for OCP($\boldsymbol{\mu}$)s}
\label{POD}
In order to build a reduced environment, \fb{two} of the major techniques \fb{that} have been exploited in the literature \fb{are} POD \cite{ballarin2015supremizer, burkardt2006pod, Chapelle2013, hesthaven2015certified} and greedy algorithm \cite{gerner2012certified, hesthaven2015certified, negri2015reduced, negri2013reduced, rozza2007stability}. We decided to rely on the first approach since the applicability of the latter \fb{requires an error estimator, which is} still not available for our nonlinear time dependent optimization problem.
\\ We now describe the POD algorithm which consists in two phases: an exploratory process based on a sample in the parameter space, in order to generate $N_{\text{max}}$ snapshots, and a compressing stage, where the snapshots are manipulated and $N < N_{\text{max}}$ basis functions are generated with the elimination of redundant information.
We provide the algorithm description as proposed in \cite{ballarin2015supremizer, burkardt2006pod, Chapelle2013, hesthaven2015certified}.
First of all, a discrete subset of parameters $\Cal P_{N_{\text{max}}} \subset \Cal P$ is chosen. If we compute truth solutions for $\boldsymbol \mu \in \Cal P_{N_{\text{max}}}$ we obtain the following sampled manifold:
$$
\mathscr M_h^{\Cal N_{\text{tot}}}=  \{ (\boldsymbol v \disc(\boldsymbol{\mu}), h\disc(\boldsymbol{\mu}), \boldsymbol u\disc(\boldsymbol{\mu}),
\boldsymbol \chi \disc (\boldsymbol{\mu}), \lambda \disc (\boldsymbol \mu))\;| \; \boldsymbol{\mu} \in \Cal P_{N_{\text{max}}}\} \subset \mathscr M^{\Cal N_{\text{tot}}}.
$$
We define $N_{\text{max}}$ as the cardinality of the set $\Cal P_{N_{\text{max}}}$ and it is clear that, when $N_{\text{max}}$ is large enough, the sampled manifold $\mathscr M_h^{\Cal N_{\text{tot}}}$ is a reliable surrogate of $\mathscr M^{\Cal N_{\text{tot}}}$.
\\ 
We decided to apply  the POD algorithm separately for all the variables involved in the system: we refer to this strategy as \emph{partitioned approach}. The final goal of the POD is to provide reduced spaces of dimension $N$ which realize the minimum of the quantities:
$$
\label{crit} 
\sqrt{\frac{1}{N_{max}}
\sum_{\boldsymbol{\mu} \in \Cal P_{N_{\text{max}}}} \underset{z_N \in {{\Cal {Y}_{v}}}_{N}}{\text{min }} \norm{\boldsymbol v\disc(\boldsymbol{\mu}) - \boldsymbol z_N}_{\Cal Y_v}^2},
$$
$$
\sqrt{\frac{1}{N_{max}}
\sum_{\boldsymbol{\mu} \in \Cal P_{N_{\text{max}}}} \underset{q_N \in {\Cal Y_{h}}_N}{\text{min}} \norm{h\disc(\boldsymbol{\mu}) - q_N}_{\Cal Y_h}},
$$
$$
\sqrt{\frac{1}{N_{max}}
\sum_{\boldsymbol{\mu} \in \Cal P_{N_{\text{max}}}} \underset{\boldsymbol \tau_N \in {\Cal {U}}_{N}}{\text{min}} \norm{\boldsymbol u\disc(\boldsymbol{\mu}) - \boldsymbol \tau_N}_{\Cal U}^2},
$$
$$
\sqrt{\frac{1}{N_{max}}
\sum_{\boldsymbol{\mu} \in \Cal P_{N_{\text{max}}}} \underset{\boldsymbol \kappa_N \in {{\Cal {Y}_v}}_{N}}{\text{min }} \norm{\boldsymbol \chi\disc(\boldsymbol{\mu}) - \boldsymbol \kappa_N}_{\Cal Y_{v}}^2},
$$
$$
\sqrt{\frac{1}{N_{max}}
\sum_{\boldsymbol{\mu} \in \Cal P_{N_{\text{max}}}} \underset{\xi_N \in {\Cal Y_h}_N}{\text{min}} \norm{\lambda\disc(\boldsymbol{\mu}) - \xi_N}_{\Cal Y_h}}.
$$
\no We now \fb{summarise} the POD-Galerkin procedure algorithm only for the velocity variable $\boldsymbol v (\boldsymbol \mu)$. In any case, the proposed approach can be identically used for the other four variables as well.
\\ Let us define a set of ordered parameters $\boldsymbol{\mu}_1, \dots, \boldsymbol{\mu}_{N_{max}}\in \Cal P_{N_{\text{max}}}$. To this parametric set, it will correspond an ordered ensemble of truth solutions, i.e. snapshots, $\boldsymbol v\disc(\boldsymbol{\mu}_1), \dots, \boldsymbol v\disc(\boldsymbol{\mu}_{N_{max}})$. We can now define the correlation matrix $\mbf C^{\boldsymbol v} \in \mbb R^{N_{max} \times N_{max}}$ of snapshots of the velocity state variable, i.e.:
$$
\mbf C_{ml}^{\boldsymbol v} = \frac{1}{N_{max}}({\boldsymbol v}\disc(\boldsymbol{\mu}_m),\boldsymbol v\disc(\boldsymbol{\mu}_l))_{\Cal {Y}}, \hspace{1cm} 1 \leq m,l \leq N_{max}.
$$
\no The next step is to solve the following eigenvalue problem
$$
\mbf C^{\boldsymbol v} x_n^{\boldsymbol v} = \theta_n^{\boldsymbol v} x_n^{\boldsymbol v}, \hspace{1cm} 1 \leq n \leq N, 
$$ 
\no with $\norm {x_n^{\boldsymbol v}}_{\Cal {Y}} = 1$. 
\fb{Assuming that the eigenvalues $\theta_1^{\boldsymbol v}, \dots, \theta_{N_{\text{max}}}^{\boldsymbol v}$ are sorted in decreasing order
we retain only the first $N$ ones, namely $\theta_1^{\boldsymbol v}, \dots, \theta_{N}^{\boldsymbol v}$, and the corresponding eigenvectors
$x_1^{\boldsymbol v}, \dots, x_{N}^{\boldsymbol v}$.}
We can now build ordered  basis functions $\{ \zeta_1^{\boldsymbol v}, \dots, \zeta_N^{\boldsymbol v}\}$ spanning the reduced space ${\Cal Y_{ v}}_N$. 
Defining  $(x_n^{\boldsymbol v})_m$ \fb{as the} \emph{m-th} component of the state velocity eigenvector $x_n^{\boldsymbol v} \in \mbb R^M$, the basis functions are given by the following relation:
$$
\zeta_n^{\boldsymbol v} = \displaystyle \frac{1}{\sqrt{{\theta_m^{\boldsymbol v}}}}\sum_{m = 1}^{N_{max}} (x_n^{\boldsymbol v})_m\boldsymbol v\disc(\boldsymbol{\mu}_m), \hspace{1cm} 1 \leq n \leq N.
$$
 Even if we performed a different POD algorithm for all the involved variables, we have not separate time instances: indeed, the snapshots are FE  solutions including all the \fb{considered temporal steps}.
This strategy is consistent with respect to the space-time FE full order discretization introduced in Section \ref{FEM}. \\
%\A{Choosing the value of $N_{max}$ is an issue related to this POD-Galerkin approach. Indeed, the accuracy of the reduced optimal solution is not guaranteed unless the number of snapshots $N_{max}$ is sufficient to reliably represent how the problem dynamics chances with respect to the parameters \cite{benner2014model, kunisch2015uniform}. However, the covariance matrix eigenvalues can guide the choice of the value $N_{max}$. Indeed, the following relation holds \cite{hesthaven2015certified,quarteroni2015reduced}:
%\begin{equation}
%\label{eq:thetas_Nmax}
%\sqrt{\frac{1}{N_{\text{max}}}
%\sum_{m = 1}^{N_{\text{max}}}  \norm{\boldsymbol v \disc(\boldsymbol{\mu}_{m}) - \mathbb P_N(\boldsymbol v\disc(\boldsymbol{\mu}_m)) }_{\mathcal Y_v}^2} = \sqrt{
%\sum_{m = N + 1}^{N_{\text{max}}}\theta_m^{\boldsymbol v},}
%\end{equation}
%where $\mathbb P_N: \mathcal Y_v \rightarrow  {\mathcal Y_v}_N$ projects functions in $\mathcal Y_v$ onto the reduced space ${\mathcal Y_v}_N$. Namely, $N_{max}$, and also $N$ can be chosen in such a way the right hand side of \eqref{eq:thetas_Nmax} is acceptably small.} 
\B{ In the proposed framework, there is no a--priori knowledge on how far the POD-besed reduced solution is from the truth approximation. Indeed, the accuracy of the reduced optimal solution is not guaranteed unless the number of snapshots $N_{max}$ is sufficient to reliably represent how the problem dynamics chances with respect to the parameters \cite{benner2014model, kunisch2015uniform}. Despite the heuristic nature behind the application of POD, the algorithm is of very common use for its great versatility since it adapts to very complex problems too, from linear to nonlinear ones, for steady or time dependent settings.}
In the next Section, we analyze how the \fb{obtained} basis functions, i.e. $\{\zeta^{\boldsymbol v}_n\}_{n=1}^N$, $\{\zeta^{h}_n\}_{n=1}^N,$ $ \{\zeta^{\boldsymbol u}_n\}_{n=1}^N,$  $ \{\zeta^{\boldsymbol \chi}_n\}_{n=1}^N$ and $\{\zeta^{\lambda}_n\}_{n=1}^N$ have to be manipulated in order to guarantee the solvability of system \eqref{ROM_optimality_system_SWE}.
\subsection{Aggregated Spaces Approach}
\label{aggr}
It is very well known in the literature \fb{that} linear PDEs constrained optimization leads to the \fb{solution} of a system in saddle point formulation \cite{Benzi, bochev2009least, hinze2008optimization, Stoll}. The saddle point framework can be extended also for lienar time dependent problems \cite{HinzeStokes, HinzeNS, Stoll1, Stoll}. In Section \ref{FEM} we showed that the saddle point framework is typical also of the linearized nonlinear system \eqref{Frechet}. The main point of solving problems based on this structure is to guarantee the inf-sup condition for the matrix $\mathsf B$, which represents the state equation. In other words, for every $\boldsymbol \mu \in \Cal P$, we want to verify the following inequality:
\begin{equation}
\label{FE_infsup}
\inf_{0 \neq \bar p} \sup_{0 \neq \bar x} \frac{\bar p^T\mathsf B \bar x}{\norm{\bar x}_{\Cal X}\norm{\bar p}_{\Cal Y}} \geq  \beta \disc(\bmu) > 0,
\end{equation}
see \cite{Babuska1971,boffi2013mixed, Brezzi74} as references. At the space-time level, relation \eqref{FE_infsup} actually holds thanks to the hypothesis on the coincidence between state and adjoint discretized spaces, which is guaranteed by the same assumption at the continuous level, as we introduced in Section \ref{ocp_SWEs}.
Now, let us suppose to have applied standard POD described in Section \ref{POD} and have obtained the following basis matrices:
$$
Z_{\bar x} = 
\begin{bmatrix}
Z_{\boldsymbol v} \\
Z_{h} \\
Z_{\boldsymbol u}
\end{bmatrix},
\spazio 
Z_{\bar p} = 
\begin{bmatrix}
Z_{\boldsymbol \chi} \\
Z_{\lambda} \\
\end{bmatrix}
\spazio \text{and} \spazio
Z =
\begin{bmatrix}
Z_{\bar x} \\
Z_{\bar p}
\end{bmatrix},
$$
where 
$
Z_{\bigstar} = [\zeta_{1}^{\bigstar} | \cdots | \zeta_{N}^{\bigstar}] \in \mathbb R^{N_t\Cal N_{\bigstar} \times N},
$ for $\bigstar = \boldsymbol v, h, \boldsymbol u, \boldsymbol \chi, \lambda$.
In order to solve the optimality system in an algebraic low dimensional framework a Galerkin projection is performed into the reduced spaces: indeed, this leads to a reduced system of the form
\begin{equation}
\label{G_compact_ROM}
\Cal G_{N}(\bar w_N; \boldsymbol \mu) \bar w_N = \bar f_N,
\end{equation}
where %\todo{ho cambiato la notazione di G}
$$\Cal G_{N}(\bar w_N; \boldsymbol \mu) := Z^T \Cal G(Z\bar w_N; \boldsymbol \mu) 
\spazio \text{and} \spazio \bar f_N := Z^T \bar f.$$
%\todo{la definizione di G N e' sbagliata. Manca sicuramente una Z nel primo argomento, e la pre e post moltiplicazione va fatta per J, non per G. Vedi todo nella sezione FE prima di sistemare}
The system \eqref{G_compact_ROM} is nonlinear, in order to solve it, we applied Newton's method and we iteratively perform
\begin{equation}
\bar {w}_N^{j + 1} := \bar {w}_N^j+ \mathbb J_N(\bar w_N^{j}; \boldsymbol \mu)^{-1}(\bar f_N - \Cal G(\bar w_N^j; \boldsymbol \mu)), \spazio j \in \mathbb N,
\end{equation}
 where the Frech\'et derivative preserves the saddle point structure, i.e. 
\begin{equation}
\label{Frechet_ROM}
\mathbb J_N (\bar w_N; \boldsymbol \mu) \bar w_N = 
\begin{bmatrix}
\mathsf A_N & \mathsf B_N^T \\
\mathsf B_N & 0 \\
\end{bmatrix}
\begin{bmatrix}
\bar x_N \\
\bar p_N
\end{bmatrix},
\end{equation}
with $J_N (\bar w_N; \boldsymbol \mu) = Z^T \mathbb J (Z \bar w_N; \bmu)Z$, $\mathsf A_N = Z_{\bar x}^T \mathsf A Z_{\bar x}$, $\mathsf B_N = Z_{\bar p}^T \mathsf B Z_{\bar x}$, $\bar x_N := Z_{\bar x}^T \bar x$ and $\bar p_N := Z_{\bar p}^T \bar p$. Also at the reduced level, equation \eqref{Frechet_ROM} implies that we still need to verify a \emph{reduced inf-sup condition} of the form
\begin{equation}
\label{ROM_infsup}
\inf_{0 \neq \bar p_N} \sup_{0 \neq \bar x_N} \frac{\bar p_N^T\mathsf B_N \bar x_N}{\norm{\bar x_N}_{\Cal X}\norm{\bar p_N}_{\Cal Y}} \geq  \beta_N(\bmu) > 0.
\end{equation}
Once again, the inequality is verified if the reduced space for state and adjoint variables is the same. \fb{However}, the introduced POD procedure will not necessarily result in the same reduced order approximation for state and adjoint. Indeed, the application of the strategy described in Section \ref{POD} will lead to the following spaces: 
\begin{align*}
& {{\Cal Y}_{\boldsymbol v}}_{N} = \text{span}\{\zeta^{\boldsymbol v}_n, \; n = 1, \dots, N\}, \\
& {{\Cal Y}_{h}}_{N} = \text{span}\{\zeta^{h}_n, \; n = 1, \dots, N\}, \hspace{2mm}\\
& {{\Cal Y}_{\boldsymbol u}}_{N} = \text{span}\{\zeta^{\boldsymbol u}_n, \; n = 1, \dots, N\}, \\
& {{\Cal Y}_{\boldsymbol \chi}}_{N} = \text{span}\{\zeta^{\boldsymbol \chi}_n, \; n = 1, \dots, N\}, \\
& {{\Cal Y}_{\lambda}}_{N} = \text{span}\{\zeta^{\lambda}_n, \; n = 1, \dots, N\}. \\
\end{align*}
Let $\Cal Q_N$ be the product space of ${{\Cal Y}_{\boldsymbol \chi}}_{N}$ and
${{\Cal Y}_{\lambda}}_{N}$: in other words, the POD defines an adjoint space $\Cal Q_N \neq \Cal Y_N$, even if the state and the adjoint space are assumed to be the same at the continuous level. It is clear that, as already specified for the continuous and discretized system, in order to guarantee the solvability of the reduced optimality system, we have to build our reduced spaces in such a way the basis functions can describe state variables as well as adjoint variables. This goal is reached thanks to the aggregated spaces technique as presented in \cite{dede2010reduced, negri2015reduced,negri2013reduced}. The main purpose of this strategy is to build a space that can be used both for state and adjoint variables. Then, let us define the aggregated spaces
$$
{{\Cal Z}_{\boldsymbol v}^{\boldsymbol \chi}}_{N}= \text{span }\{\zeta^{\boldsymbol v}_n, \zeta^{\boldsymbol \chi}_n, \; n = 1, \dots, N\}
\hspace{2mm} \text{and} \hspace{2mm} 
{{\Cal Z}_{h}^{\lambda}}_{N}= \text{span }\{\zeta^{h}_n, \zeta^{\lambda}_n, \; n = 1, \dots, N\}.
$$
The product space $\Cal Z_N = {{\Cal Z}_{\boldsymbol v}^{\boldsymbol \chi}}_{N} \times {{\Cal Z}_{h}^{\lambda}}_{N}$ can actually give a representation of the reduced state variable $(\boldsymbol v_N(\boldsymbol \mu), h_N(\bmu))$ and the reduced adjoint variable 
$(\boldsymbol \chi_N(\boldsymbol \mu), \lambda_N(\bmu))$. Moreover, setting $\Cal Y_N \equiv \Cal Q_N \equiv \Cal Z_N$, the inf-sup condition \eqref{ROM_infsup} holds.
Concerning the control \fb{function} space, a standard POD-procedure can be applied, building
$$
\Cal U_N = \text{span }\{ \zeta^{\boldsymbol u}_n, \; n = 1, \dots, N\}.
$$
The aggregated space technique allows us to define new basis matrices of the form: \\
$
Z_{\boldsymbol v} \equiv Z_{\boldsymbol \chi} = [\zeta_{1}^{ \boldsymbol v} | \cdots | \zeta_{N}^{ \boldsymbol v}| \zeta_{1}^{ \boldsymbol \chi} | \cdots | \zeta_{N}^{ \boldsymbol \chi}] \in \mathbb R^{N_t\Cal N_{\boldsymbol v} \times 2N},
$
$
Z_{\boldsymbol u} = [\zeta_{1}^{ \boldsymbol u} | \cdots | \zeta_{N}^{ \boldsymbol u}] \in \mathbb R^{N_t\Cal N_{\boldsymbol u} \times N}
$
and 
$
Z_{h} \equiv Z_{\lambda} = [\zeta_{1}^{h} | \cdots | \zeta_{N}^{h}| \zeta_{1}^{\lambda} | \cdots | \zeta_{N}^{\lambda}] \in \mathbb R^{N_t\Cal N_{h} \times 2N}.
$
The new spaces are actually increasing the dimension of the system since the global reduced dimension is $N_{\text{tot}} = 9N$. Although, the strategy guarantees the reduced inf-sup condition \eqref{ROM_infsup} and, consequently, the existence of an optimal solution. \fb{Still}, 
$N_{\text{tot}} < \Cal N_{\text{tot}}$, i.e. we still work in a reduced dimensional framework. \\
\no We introduced all the notions needed in order to reduce nonlinear time dependent \ocp s. Anyway, we still miss fundamental assumptions which allow ROM to be very advantageous under the point of view of computational costs: it will be the topic of the next Section.
%%%%%%%%%%%%%%%%%%%%%%%%%%%
%%%%%%%%%%%%%%%%%%%%%%%%%%%%%%
%%%%%%%%%%%%%%%%%%%%%%%%%
%

\subsection{Efficient ROM and Affinity Assumption: Offline--Online decomposition}
\label{aff}
Exploiting a reduced strategy is convenient only if fast simulations can be assured in order to analyse different configurations of the physical system for several parameters. To guarantee an efficient applicability of ROM techniques, the system is assumed to be affinely decomposed. In other words, all the quantities involved in the system have to be interpreted as the product of $\boldsymbol \mu -$ dependent quantities and $\boldsymbol \mu -$independent quantities, i.e.
the equations involved can be recast as:
\begin{equation*}
D_{\boldsymbol v}\Lg((\boldsymbol v, h), \boldsymbol u, (\boldsymbol \chi, \lambda))[\boldsymbol z] =
\displaystyle  \sum_{q=1}^{Q_{D_{\boldsymbol v}\Lg}} \Theta_{D_{\boldsymbol v}\Lg}^q(\boldsymbol{\mu})D_{\boldsymbol v}\Lg^q((\boldsymbol v, h), \boldsymbol u, (\boldsymbol \chi, \lambda))[\boldsymbol z],
\end{equation*}
\begin{equation*}
D_{h}\Lg((\boldsymbol v, h), \boldsymbol u, (\boldsymbol \chi, \lambda))[q] =
\displaystyle  \sum_{q=1}^{Q_{D_{h}\Lg}} \Theta_{D_{h}\Lg}^q(\boldsymbol{\mu})D_{h}\Lg^q((\boldsymbol v, h), \boldsymbol u, (\boldsymbol \chi, \lambda))[q],
\end{equation*}
\begin{equation}
\label{affinity}
D_{\boldsymbol u}\Lg((\boldsymbol v, h), \boldsymbol u, (\boldsymbol \chi, \lambda))[\boldsymbol \tau] =
\displaystyle  \sum_{q=1}^{Q_{D_{\boldsymbol u}\Lg}} \Theta_{D_{\boldsymbol u}\Lg}^q(\boldsymbol{\mu})D_{\boldsymbol u}\Lg^q((\boldsymbol v, h), \boldsymbol u, (\boldsymbol \chi, \lambda))[\boldsymbol \tau],
\end{equation}
\begin{equation*}
D_{\boldsymbol \chi}\Lg((\boldsymbol v, h), \boldsymbol u, (\boldsymbol \chi, \lambda))[\boldsymbol \kappa] =
\displaystyle  \sum_{q=1}^{Q_{D_{\boldsymbol \chi}\Lg}} \Theta_{D_{\boldsymbol \chi}\Lg}^q(\boldsymbol{\mu})D_{\boldsymbol \chi}\Lg^q((\boldsymbol v, h), \boldsymbol u, (\boldsymbol \chi, \lambda))[\boldsymbol \kappa],
\end{equation*}
\begin{equation*}
D_{\lambda}\Lg((\boldsymbol v, h), \boldsymbol u, (\boldsymbol \chi, \lambda))[\xi] =
\displaystyle  \sum_{q=1}^{Q_{D_{\lambda}\Lg}} \Theta_{D_{\lambda}\Lg}^q(\boldsymbol{\mu})D_{\lambda}\Lg^q((\boldsymbol v, h), \boldsymbol u, (\boldsymbol \chi, \lambda))[\xi].
\end{equation*}
\no  for some finite $Q_{D_{\boldsymbol v}\Lg}, Q_{D_{h}\Lg}, Q_{D_{\boldsymbol u}\Lg},$ $ Q_{D_{\boldsymbol \chi}\Lg}, Q_{D_{\lambda}\Lg}$, where $\Theta_{D_{\boldsymbol v}\Lg}^q, \Theta_{D_{h}\Lg}^q, \Theta_{D_{\boldsymbol u}\Lg}^q, $ $\Theta_{D_{\boldsymbol \chi}\Lg}^q$ and $\Theta_{D_{\lambda}\Lg}^q$ are $\boldsymbol{\mu}-$dependent smooth functions, whereas $D_{\boldsymbol v}\Lg^q((\boldsymbol v, h), \boldsymbol u, (\boldsymbol \chi, \lambda)),$ \\$ D_{h}\Lg^q((\boldsymbol v, h), \boldsymbol u, (\boldsymbol \chi, \lambda)), $ $D_{\boldsymbol u}\Lg^q((\boldsymbol v, h), \boldsymbol u, (\boldsymbol \chi, \lambda)),$ $ D_{\boldsymbol \chi}\Lg^q((\boldsymbol v, h), \boldsymbol u, (\boldsymbol \chi, \lambda))$ and \\$D_{\lambda}\Lg^q((\boldsymbol v, h), \boldsymbol u, (\boldsymbol \chi, \lambda))$ are $\boldsymbol{\mu} -$independent quantities describing the optimality system. Thanks to affine decomposition, the \fb{solution} of an OCP($\bmu$) can be performed in two different steps: an \textbf{offline} stage which consists in assembling all the parameter independent quantities, building the reduced \fb{function} spaces and storing all the mentioned quantities; then, an \textbf{online} stage deals with all the $\boldsymbol{\mu}-$dependent quantities and solves the whole reduced system \eqref{ROM_optimality_system_SWE}. The latter phase is performed at every new parameter evaluation and gives us information about physical configurations in a small amount of time since it does not depend on the discrete full order dimension $\Cal N_{\text{tot}}$. 
%Although, the drawback could be the expensive computational effort required for the offline phase: anyway, it is still performed only once and the storing process guarantees an efficient real-time computation in a many query context.
We underline that the case of interest deals at most with only quadratically
nonlinear terms, we can guarantee the affinity assumption storing the appropriate
nonlinear terms in third order tensors.
\no If the OCP($\bmu$) does not fulfill the decomposition \eqref{affinity}, the empirical interpolation method (EIM) can recover the assumption, see \cite{barrault2004empirical} or \cite[Chapter 5]{hesthaven2015certified} as references. \\
In the next Section, we are going to present some numerical results for a nonlinear time-dependent OCP($\bmu$) governed by SWE equations, in order to assert the applicability of ROM for this complicated model, which can be of great interest for in many fields of natural sciences and engineering.

\section{Numerical Results}
\label{Results}
This Section aims at validating the numerical performances of POD-Galerkin projection over a nonlinear time dependent OCP($\bmu$) governed by SWEs. We put in a parametrized optimal control framework the academic test case presented in \cite{ferrari2004new}. The experiment can be read as an inverse problem on the forcing term, i.e. find the control variable, needed in order to have a desired velocity-height state profile. Given a parameter 
$\bmu =(\mu_1, \mu_2, \mu_3)$ in the parametric space $\Cal P = (0.00001,1.) \times (0.01, 0.5) \times (0.1, 1.)$, we solved the optimality system \eqref{optimality_system_SWEs} built through  Lagrangian approach, in the fashion of optimize-then-discretize technique, over the water basin described by $\Omega = (0, 10) \times (0,10)$. As we already discussed in Section \ref{SWE}, the flat bathymetry does not actually affect the system, so we used $z_b = 0$. We simulate our system evolution in the time interval $(0,T)$ with T = 0.8s. The OCP($\bmu$) simulates the spreading of a mass of water with an initial Gaussian distributed elevation and null initial velocity: i.e.
$$
\boldsymbol v_0 = \boldsymbol 0, \spazio \text{and} \spazio h_0 = 0.2(1 + 5e^{(-(x_1 - 5)^2 - (x_2 - 5)^2 + 1))}),
$$
where $x_1$ and $x_2$ are the spatial coordinates. Under a controlling forcing term representing wind action and bottom friction, we want our solution to be similar to $(\mu_3 \boldsymbol v_d, \mu_3 h_d)$, where  $(\boldsymbol v_d, h_d)$ is the solution of the uncontrolled state equation \eqref{SWEs_state}, with null initial velocity and initial elevation 
${h_d}_{0} = 2e^{(-(x_1 - 5)^2 - (x_2 - 5)^2 + 1)}$
and null forcing term $\boldsymbol u = \boldsymbol 0$. \\
\no In Section \ref{SWE}, we already specified the diffusive and advective role of $\mu_1$ and $\mu_2$.
In Table \ref{table_data} we report all the specifics of the experiment that we are going to describe.
\no The goal of the presented optimal control problem is to make our solution $(\boldsymbol v, h)$ the most similar to the desired above mentioned profile. Once again, we underline that we work in a parametrized framework, i.e. the controlled solution changes for different values of $\bmu \in \Cal P$. All the results we present are given by the parameter 
$\bmu = (0.1, 0.5, 1.)$. Following the space-time discretization technique proposed in \cite{saleri2007geometric},  we used linear polynomials for the truth approximation of all the variables, i.e. $r_{\boldsymbol v} = r_{h} = r_{\boldsymbol u} = 1$. With respect to time discretization, we divided the time interval with  $\Delta t = 0.1$, which leads to a number of time steps $N_t = 8$. The problem solved is quite complex even with this small amount of time steps. In any case, $\Delta t$ can be reduced following the iterative techniques exploited in \cite{HinzeStokes, HinzeNS, Stoll1, Stoll}. Although, for the sake of simplicity, we exploited a direct solver for the algebraic system \eqref{G_compact}. In the end, at the truth approximation level, we deal with a system of a total dimension $\Cal N_{\text{tot}} = N_t \times {\Cal {N}} = 94'016$. 

\begin{figure}[H]
\begin{center}
\includegraphics[scale = 0.22]{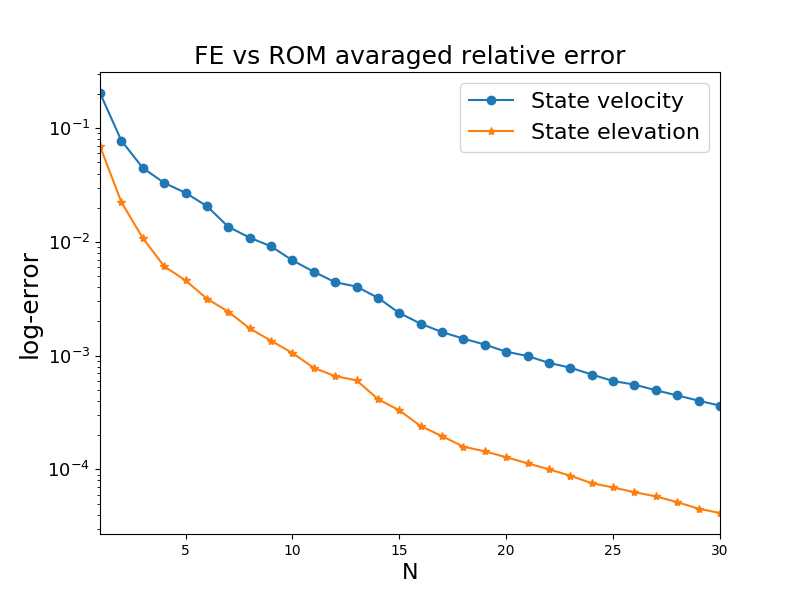}
\includegraphics[scale = 0.22]{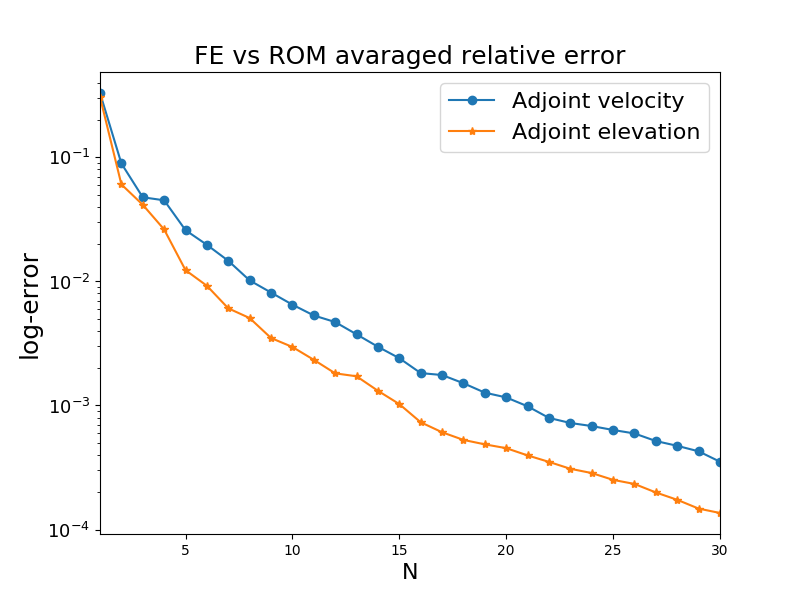}
\includegraphics[scale = 0.22]{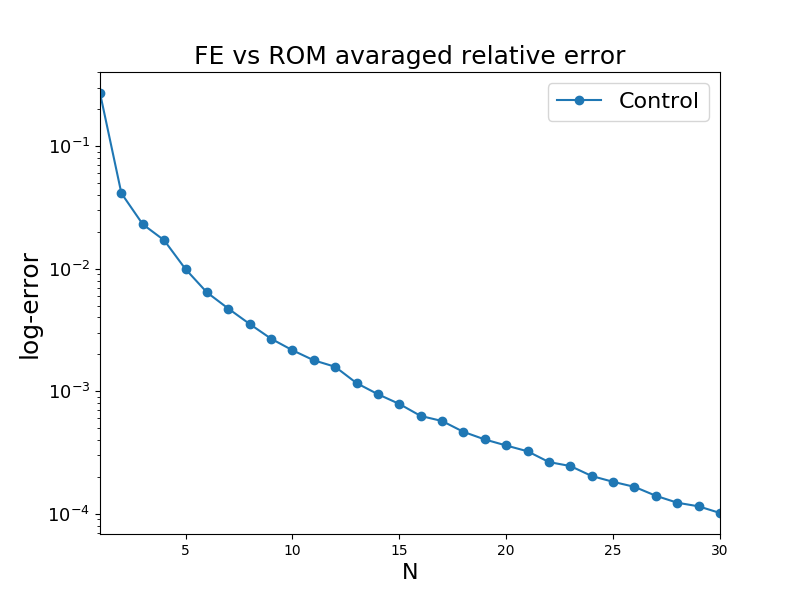}
\hspace{-1cm}\caption{Averaged relative error between space-time and ROM approximation for state velocity and elevation profile (top left), adjoint velocity and elevation profile (top right), and control (bottom).}
\label{errors}
\end{center}
\end{figure}

\begin{figure}[H]
\begin{center}
\includegraphics[scale = 0.12]{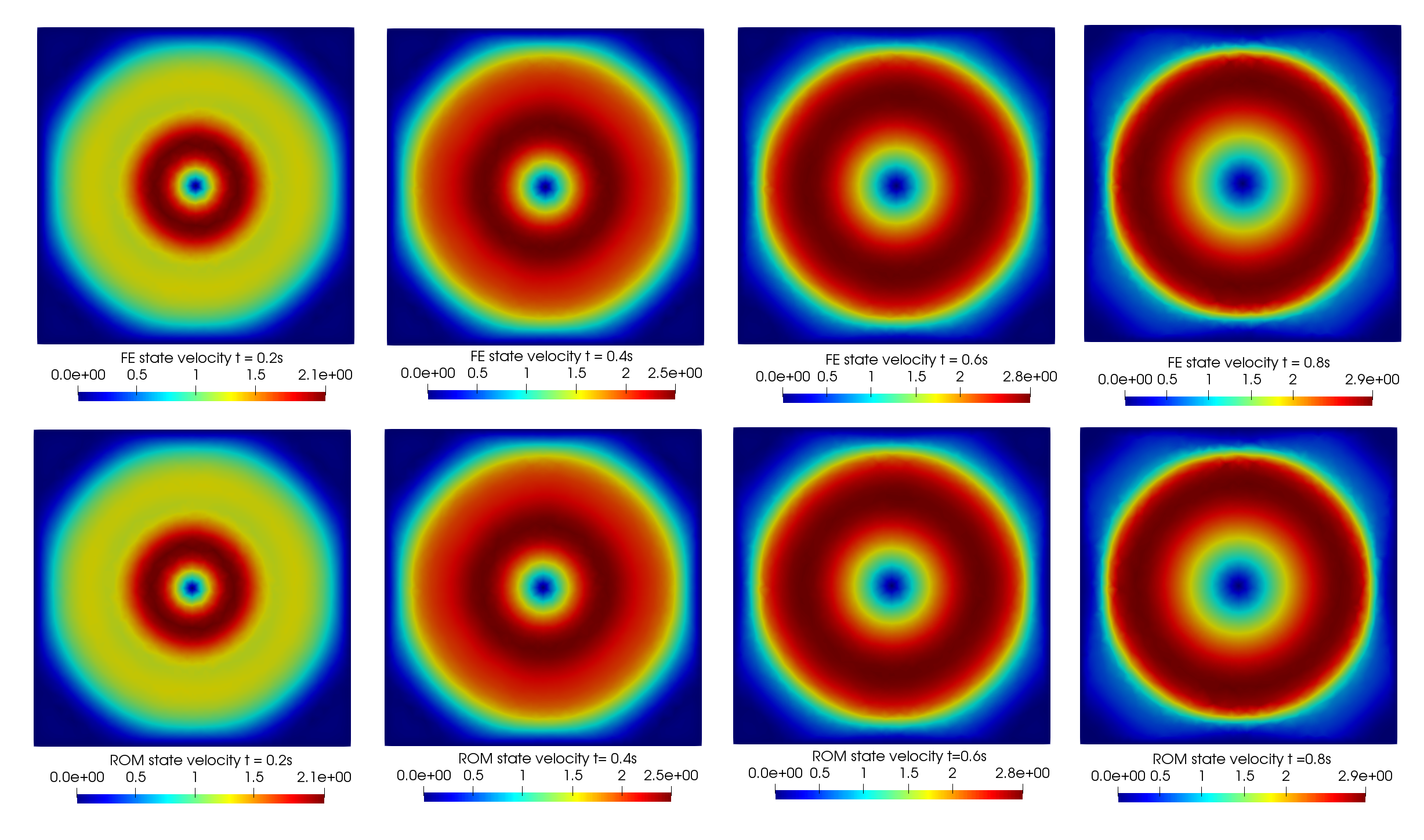}
\hspace{-1cm}\caption{FE  state velocity profile (top) compared to ROM state velocity profile (bottom) for $t = 0.2s, 0.4s, 0.6s, 0.8s$ and for $\bmu = (0.1, 0.5, 1)$.}
\label{v_comp}
\end{center}
\end{figure}

\begin{figure}[H]
\begin{center}
\includegraphics[scale = 0.12]{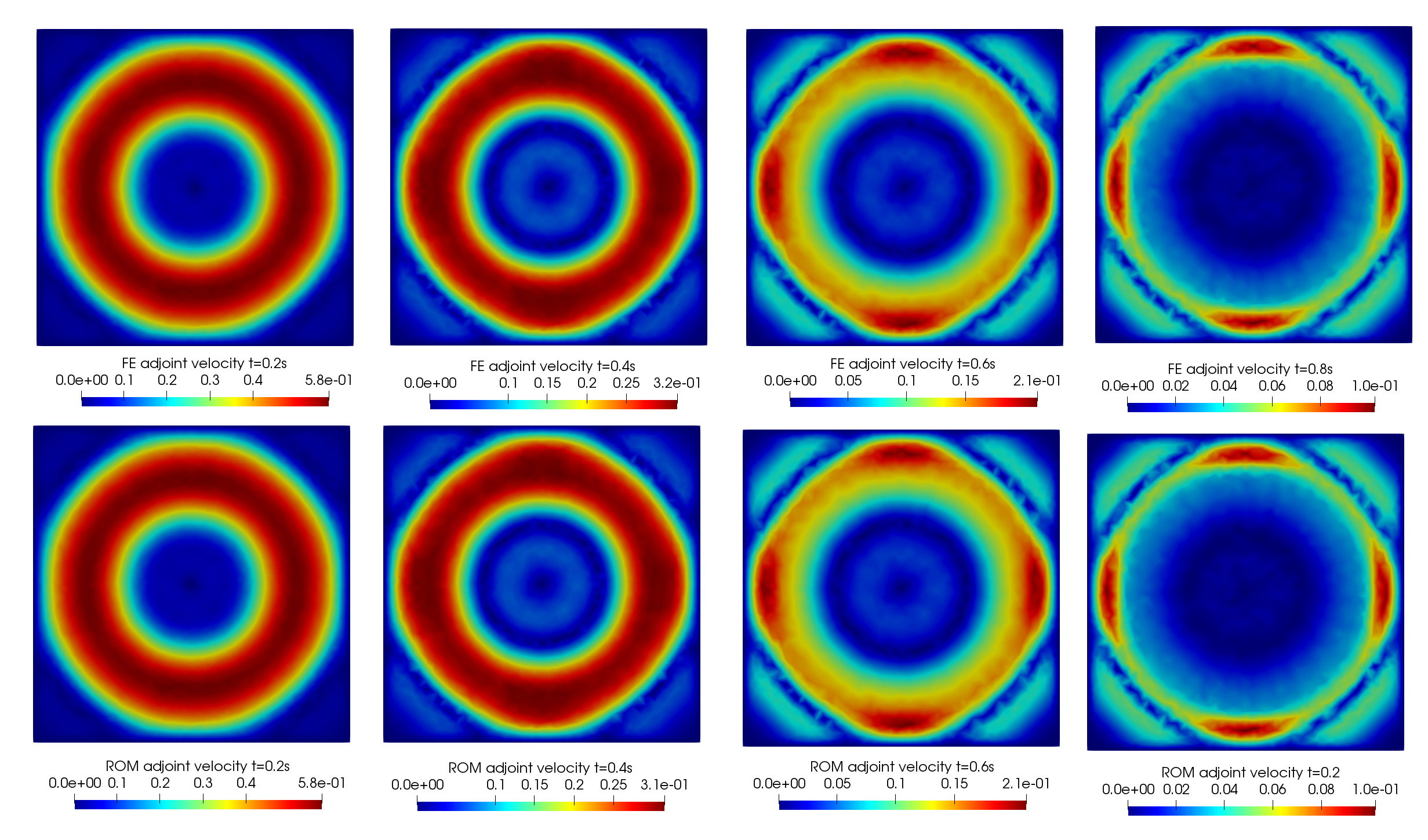}
\hspace{-1cm}\caption{FE  adjoint velocity profile (top) compared to ROM adjoint velocity profile (bottom) for $t = 0.2s, 0.4s, 0.6s, 0.8s$ and for $\bmu = (0.1, 0.5, 1)$.}
\label{w_comp}
\end{center}
\end{figure}

\begin{figure}[H]
\begin{center}
\includegraphics[scale = 0.12]{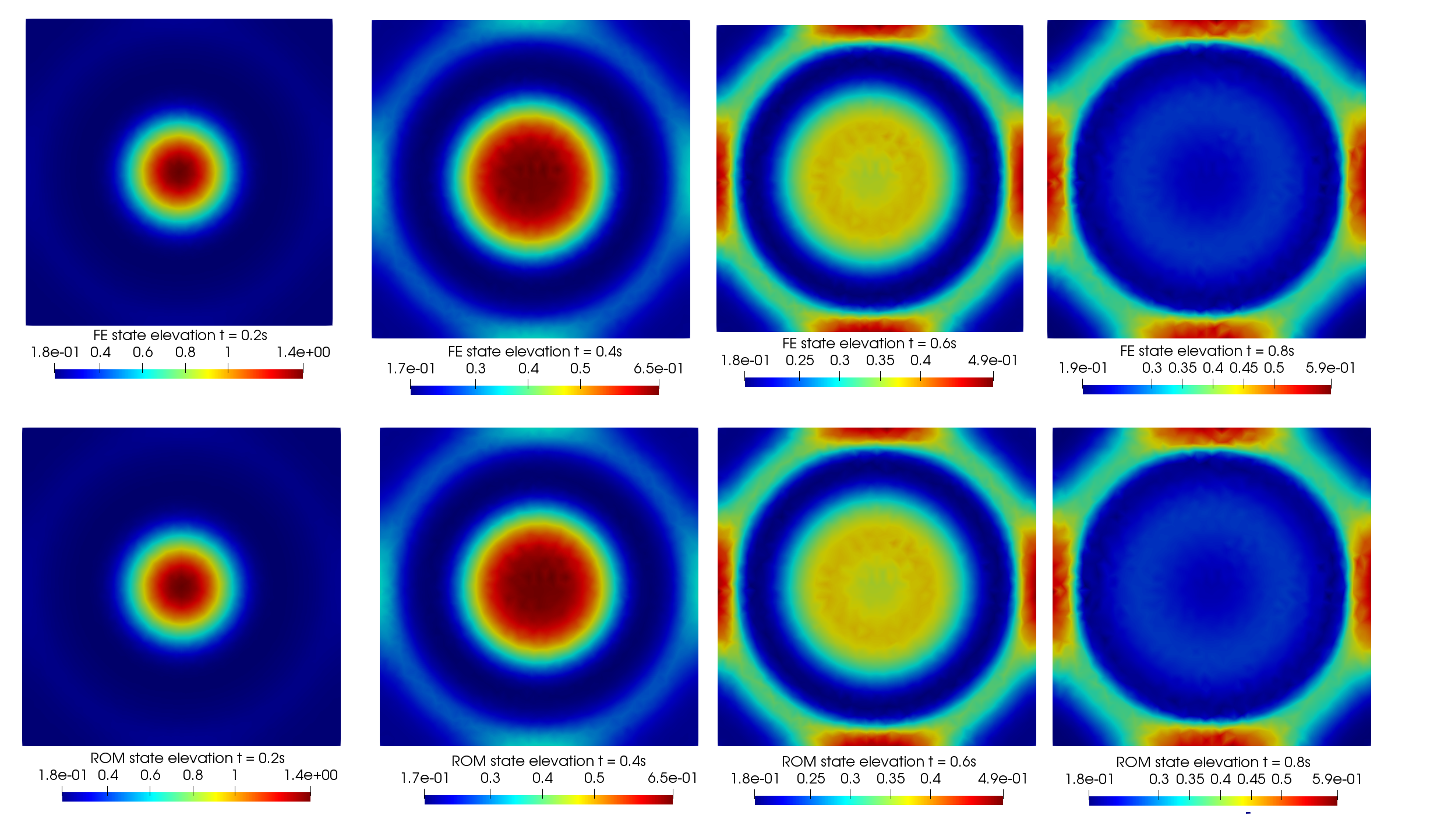}
\hspace{-1cm}\caption{FE  state elevation profile (top) compared to ROM state elevation profile (bottom) for $t = 0.2s, 0.4s, 0.6s, 0.8s$ and for $\bmu = (0.1, 0.5, 1)$.}
\label{h_comp}
\end{center}
\end{figure}

\begin{figure}[H]
\begin{center}
\includegraphics[scale = 0.12]{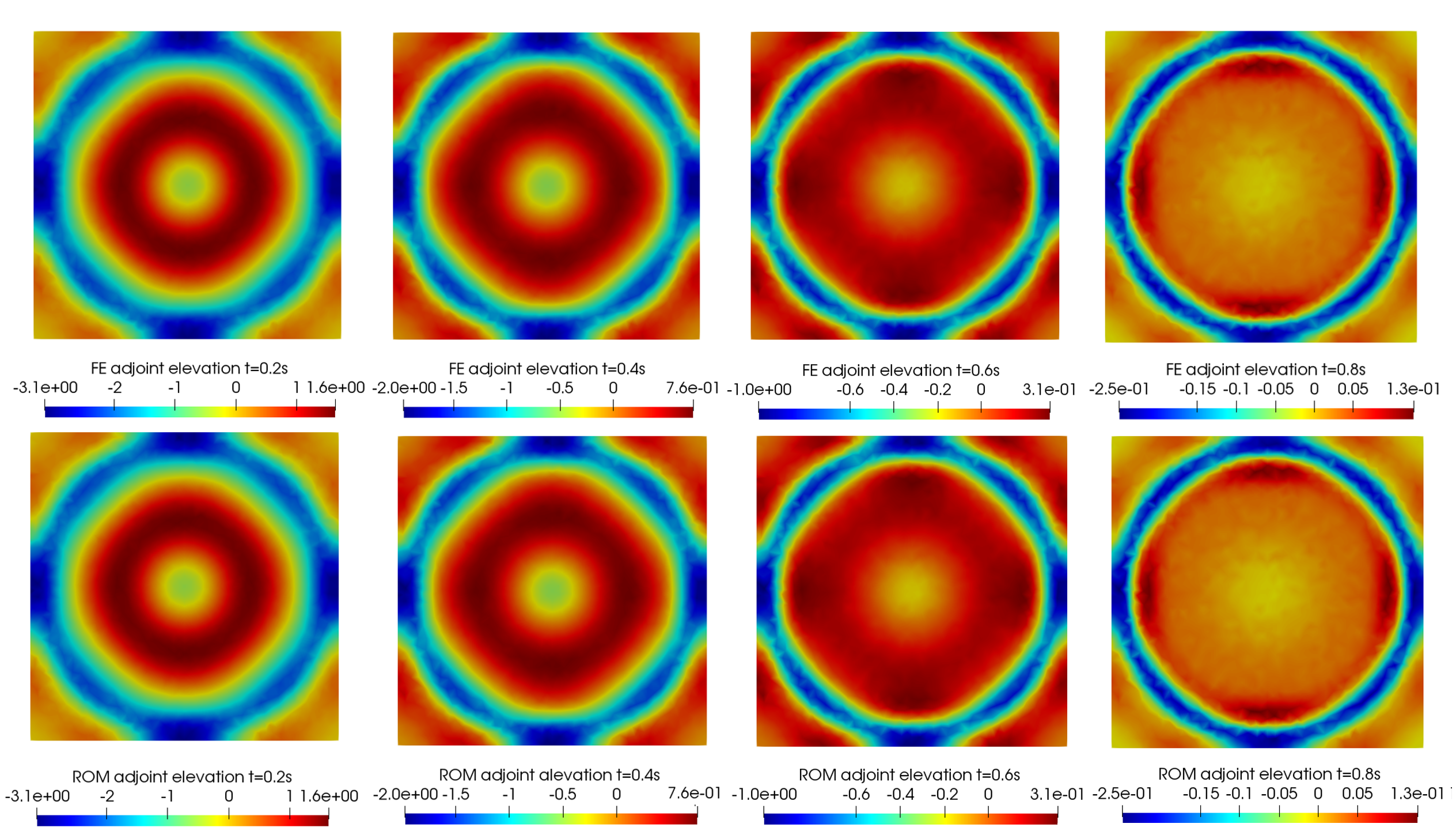}
\hspace{-1cm}\caption{FE  adjoint elevation profile (top) compared to ROM adjoint elevation profile (bottom) for $t = 0.2s, 0.4s, 0.6s, 0.8s$ and for $\bmu = (0.1, 0.5, 1)$.}
\label{q_comp}
\end{center}
\end{figure}

In order to build the reduced optimality system \eqref{ROM_optimality_system_SWE}, we applied the partitioned POD-Galerkin approach presented in Section \ref{POD}. First of all, we built five correlation matrices with $N_{\text{max}} = 100$ for all the variables, respectively. \B{The choice of $N_{\text{max}}$ is affected by an increasing effort in solving the offline phase. Indeed, the described problem has huge computational limitations both in time and in storage memory exploited for the basis construction: they drastically grow for large values of $N_{\text{max}}$. Although, the value we used led to a feasable offline phase in terms of computational time}. \\ Let us define the \emph{basis number} $N$, i.e. the number of eigenvalues/eigenvectors retained from the correlation matrices compression process of the POD. For this test case, the basis functions  were obtained choosing $N = 30$. The basis number considered allowed us to well describe the full order approximated system in the reduced framework, as the reader can notice from the average relative errors represented in Figure \ref{errors} with the following norms: 
$$
\intTime{\norm{\boldsymbol v \disc - \boldsymbol v_N}^2_{H^1}}, \hspace{2mm}
\hspace{2mm}
\intTime{\norm{h \disc - h_N}^2_{L^2}} \hspace{2mm}
$$
$$
\intTime{\norm{\boldsymbol u \disc - \boldsymbol u_N}^2_{L^2}}, \hspace{2mm}
\hspace{2mm}
\intTime{\norm{\boldsymbol \chi \disc - \boldsymbol \chi_N}^2_{H^1}}, \hspace{2mm}
\text{and }\hspace{2mm}
\intTime{\norm{\lambda \disc - \lambda_N}^2_{L^2}}.\hspace{2mm}
$$
\begin{table}[H]
\centering
\caption{Data for the \ocp  $\:$ governed by SWEs.}
\label{table_data}
\begin{tabular}{ c | c }
\toprule
\textbf{Data} & \textbf{Values} \\
\midrule
$\mathcal P$ &  $(0.00001,1) \times (0.01, 0.5) \times (0.1, 1)$\\
\midrule
$[0,T]$ &  $[0s, 0.8s]$\\
\midrule
 values of $(\mu_1, \mu_2, \mu_3,  \alpha)$ &  $(0.1, 0.5, 1, 10^{-1})$ \\
\midrule
$N_{\text{max}}$ & 100\\ 
\midrule
$N$ & 30 \\
\midrule
Sampling Distribution & Uniform \\
\midrule
$\Cal N_{\text{tot}}$ & 94'016 \\
\midrule
ROM System Dimension & 270 \\
\bottomrule
\end{tabular}
\end{table}
The errors is averaged over a testing set of $20$ parameters uniformly distributed: as expected, it decreases with respect the basis number $N$ reaching a minimum value of $10^{-3}$ state and adjoint velocity and a value of $10^{-4}$ for state and adjoint elevation together with the control variable.
%\begin{figure}[H]
%\begin{center}
%\includegraphics[scale = 0.12]{../Images/vh_uncontrolled}
%\caption{FE  uncontrolled velocity (top) and elevation (bottom) for $t = 0.2s, 0.4s, 0.6s, 0.8s$ and for $\mu_1 = 0.1$ and $\mu_2 = 0.5$.}
%\label{uncontrolled}
%\end{center}
%\end{figure}
%\begin{figure}[H]
%\begin{center}
%\includegraphics[scale = 0.12]{../Images/vh_desired}
%\hspace{-1cm}\caption{FE  desired velocity (top) and elevation (bottom) for $t = 0.2s, 0.4s, 0.6s, 0.8s$ and for $\mu_1 = 0.0001$ and $\mu_2 = 1$.}
%\label{desired}
%\end{center}
%\end{figure}
\no Moreover, the effectiveness of the reduced model can be understood also from the comparison between the space-time solutions and the ROM solutions presented in Figures \ref{v_comp} and 
\ref{h_comp} for the velocity and elevation state at $t = 0.2s, 0.4s, 0.6s, 0.8s$, respectively. The ROM procedure leads to a good representation of the space-time solutions for the different time instances considered.  
The same conclusions can be drawn for the adjoint variables for velocity and height in Figure \ref{w_comp} and \ref{q_comp}, respectively. For the sake of brevity we do not show the control variable. Indeed, the adjoint variable has the same behaviour of the control variable scaled by the factor $\alpha$, as a consequence of the optimality equation \eqref{opt_eq_k}. \no Let us analyse the computational time comparison between space-time and ROM simulations. We refer to \emph{the speedup index}: it represents how many reduced simulations can be performed in the time of one space-time optimality system \fb{solve}. \fb{The speedup depends very mildly on the value of $N$, and} it is of the order of $O(30)$ for $N = 1, \dots, 30$.\\
We remind that in order to guarantee the solvability of the reduced saddle point problem arising from the linearized system, we used aggregated space technique presented in Section \ref{aggr}: it increased the reduced dimension of the system to $N_{\text{tot}}= 9 N = 270$. Anyway, the speedup index underlines that it is actually convenient to perform a projection even in this larger reduced space, since the whole optimality system is actually very complex to be solved at the space-time level, most of all if many simulations are required in order to better study several parametric configurations. 
The next Section is dedicated to some comments and perspectives on improvements and future research focus with respect to the presented topic.

%
%\begin{figure}[H]
%\begin{center}
%\includegraphics[scale = 0.12]{../Images/u_comp}
%\hspace{-1cm}\caption{FE  control variable (top) compared to ROM control variable (bottom) for $t = 0.2s, 0.4s, 0.6s, 0.8s$ and for $\bmu = (0.1, 0.5, 1)$.}
%\label{u_comp}
%\end{center}
%\end{figure}

%The proposed example shows how reduced formulation can be effective for very complex systems in environmental marine sciences such as \ocp s governed by SWEs state equation. However, it has still some limits in its applications. Anyway, this work is the first attempt to reduce a physical parametrization of a nonlinear time dependent \ocp $\;$ governed by SWEs, at the best of our knowledge.\\

\section{Conclusion and Perspectives}

\label{conc}
In this work, we propose ROMs as a suitable tool to solve a parametrized nonlinear time dependent OCP($\bmu$) governed by SWE, a very important model widely spread in several environmental applications such as marine ecosystem management and coastal engineering. To the best of our knowledge, it is the first time that parametrized reduction is exploited for this kind of solution tracking coastal model: this work aims at showing how ROMs could be very effective in the marine management field deeply characterized by a growing demanding computational effort. Working in a low dimensional framework allows us to perform accurate simulations in a small amount of time compared to the space-time approximation. Moreover, the proposed methodology, based on a POD-Galerkin projection of the Lagrangian based optimality system, is general and can be easily applied to other nonlinear time dependent state equations. Indeed, a complete reduced data assimilated nonlinear and time dependent framework is presented and, since it can be used for general state equations, it can be exploited in several environmental problem configurations. \A{ We remark that the proposed formulation has the strength of sharing many similarities with steady problems. Besides this property, we reach the goal of showing how ROMs can be effective in this context and this was a first step towards a deeper investigation of the topic.} \\
Indeed, some possible advances to this work follow. First of all, a more theoretical analysis of the considered problem is still partially missing due to the great complexity of the state equation itself. The analysis of the parametrized optimal control framework governed by SWEs could be a topic of future investigation, as well as the development of an error estimator which could allow us to rely on more efficient ROM greedy-based algorithms. \\
\A{Another point of interest concerns the hyperbolic nature of the system. Indeed, it is well known that SWE might lead to a poor ROM approximation \cite{GRIMBERG, TaddeiSWE}, however, the example we proposed seemed not to suffer from this issue. In order to take into consideration more complex problems, a deeper analysis of the role of optimal control in this setting is needed.}\\
Moreover, for environmental sciences applications, another important development for nonlinear  time dependent problems could be the use of random input parameters as an extension of \cite{CarereStrazzulloRozza}, since, in marine ecosystems, it is not always possible to assign deterministic values for the parameters describing the physical model.

%\\ Since nonlinear time dependent \ocp s are of great importance in environmental applications, we propose as a further step the simulation of the data assimilation SWE model in an actual geographical domain, taking into consideration also a non-constant bathymetry effect. This will lead to a complete \ocp, which could give more realistic coastal predictions.
%\\ We would also like to move towards a reduced uncertainty quantification modelling for environmental \ocp s, since, in these marine ecosystems, it is not always possible to assign deterministic values for the parameters involved in the physical model.
%\\ We aim at improving the proposed reduced model \ocp in order to apply it in many-query contexts, where several simulations have to be analysed in order to study marine physical phenomena, and in real-time applications, which require fast and reliable parametrized simulations.

\appendix
\A{
\section{Details on the Continuous Model}
\label{1}
We here propose a more detailed discussion on the problem formulation introduced in Section \ref{ocp_SWEs}.
First of all,  the state equation \eqref{SWEs_state} can be interpret in weak formulation:
given a parameter $\boldsymbol \mu \in \Cal P$, find 
$x\in \mathcal X$ which minimizes \eqref{J} constrained to 
$\Cal S(x, (\boldsymbol \kappa, \xi); \boldsymbol \mu)=0$ for all $(\boldsymbol \kappa, \xi) \in \Cal Y$, where $\Cal S(x, (\boldsymbol \kappa, \xi); \boldsymbol \mu)=0$ reads as follows:
\begin{align*}
\begin{cases}
\displaystyle \intTimeSpace{\dt{\boldsymbol v} \fb{\cdot} \boldsymbol \kappa }+ \mu_1 a_1((\boldsymbol v, h), (\boldsymbol \kappa, \xi)) +
\mu_2 \displaystyle \intTime{a_1^{\text{nl}}((\boldsymbol v, h), (\boldsymbol v, h), (\boldsymbol \kappa, \xi))} 
\\ \qquad \qquad \qquad \qquad + \displaystyle \intTime{a_2((\boldsymbol v, h), (\boldsymbol \kappa, \xi)} = \displaystyle \intTimeSpace {\boldsymbol u \fb{\cdot} \boldsymbol \kappa} \\
\displaystyle \intTimeSpace{\dt{h} \xi} + \intTime{a_2^{\text{nl}}((\boldsymbol v, h), (\boldsymbol v, h), (\boldsymbol \kappa,\xi))} = 0.
\end{cases}
\end{align*}
The forms $a_1\cd, a_2\cd, a_1^{nl}(\cdot, \cdot, \cdot), a_2^{nl} (\cdot, \cdot, \cdot)$ are defined as follows:
\begin{align*}
 a_1 \goesto{\Cal Y}{\Cal Y}{\mathbb R} & \spazio a_1((\boldsymbol v, h), (\boldsymbol \kappa, \xi)) =   \displaystyle \intSpace {\nabla \boldsymbol v \fb{:} \nabla \boldsymbol \kappa},\\
a_2 \goesto{\Cal Y}{\Cal Y}{\mathbb R} & \spazio a_2((\boldsymbol v, h), (\boldsymbol \kappa, \xi)) = \displaystyle  \intSpace {g \boldsymbol \kappa \fb{:} \nabla h}, \\
a_1^{\text{nl}} : {\Cal Y} \times {\Cal Y} \times {\Cal Y} \rightarrow {\mathbb R} & \spazio  a_1^{\text{nl}}((\boldsymbol v, h), (\boldsymbol w, \varphi), (\boldsymbol \kappa, \xi)) = \displaystyle \intSpace {\boldsymbol \kappa \fb{:} ( \boldsymbol v \cdot \nabla)\boldsymbol w},\\
a_2^{\text{nl}} : {\Cal Y} \times {\Cal Y} \times {\Cal Y} \rightarrow {\mathbb R} & \spazio  a_2^{\text{nl}}((\boldsymbol v, h), (\boldsymbol w, \varphi), (\boldsymbol \kappa, \xi)) = \displaystyle  \intSpace {\xi \dive(\varphi \boldsymbol v)}.\\
\end{align*}
The nonlinear nature of the bilinear forms have been specified through the apex ``\emph{nl}''.
Considering the adjoint variable $(\boldsymbol \chi, \lambda) := (\boldsymbol \chi(\bmu), \lambda(\bmu)) \in \Cal Y$, it is clear that the weak form of the optimality equation is given by  
$$\alpha \intTimeSpace {\boldsymbol u \fb{\cdot} \boldsymbol \tau} = \intTimeSpace {\boldsymbol \chi \fb{\cdot} \boldsymbol \tau} \spazio \forall \tau \in \Cal U,$$
while the weak {adjoint equation}, reads:
\begin{equation*}
\label{adj}
\begin{aligned}
\begin{cases}
\displaystyle 
\intTimeSpace{\boldsymbol v \fb{\cdot} \boldsymbol z} - \intTimeSpace{\dt{\boldsymbol \chi} \fb{\cdot} \boldsymbol z }+ \mu_1 \intTime{a_1^{\ast}((\boldsymbol \chi, \lambda), (\boldsymbol z, q))} \\ \displaystyle \qquad  \qquad \qquad +
\mu_2 \intTime{{a_1^{\text{nl}}}^{\ast}((\boldsymbol v, h), (\boldsymbol \chi, \lambda), (\boldsymbol z, q))} \\
\displaystyle \qquad \qquad \qquad \qquad \qquad+ \intTime{{a_2^{\text{nl}}}\dual((\boldsymbol v, h), (\boldsymbol \chi, \lambda), (\boldsymbol z, q))}= \intTimeSpace{\boldsymbol v_d \fb{\cdot} \boldsymbol z}, \vspace{1mm}\\
\displaystyle \intTimeSpace{h q}  - \intTimeSpace{\dt{h} q} + \intTime{{a_3^{\text{nl}}}^{\ast}((\boldsymbol v, h), (\boldsymbol \chi, h), (\boldsymbol z, q))}  \\ \qquad \qquad \qquad \hspace{3cm} \displaystyle +\intTime{ a_3^{\ast}((\boldsymbol \chi, \lambda), (\boldsymbol z, q))}= \intTimeSpace{h_d q}, 
\end{cases}
\end{aligned}
\end{equation*}
for all $(\boldsymbol z, q) \in \Cal Y$, where the involved forms are defined as
\begin{align*}
 a_1\dual \equiv a_1 \goesto{\Cal Y}{\Cal Y}{\mathbb R} & \qquad a_1((\boldsymbol \chi, \lambda), (\boldsymbol z, q)) =   \intSpace {\nabla \boldsymbol \chi \fb{:} \nabla \boldsymbol z},\\
{a_1^{\text{nl}}}^{\ast} : {\Cal Y} \times {\Cal Y} \times {\Cal Y} \rightarrow {\mathbb R} & 
\qquad  {a_1^{\text{nl}}}^{\ast}((\boldsymbol v, h), (\boldsymbol \chi, \lambda), (\boldsymbol z, q)) =  - \intSpace{(\boldsymbol v \cdot \nabla)\boldsymbol \chi \fb{\cdot} \boldsymbol z} \\ 
& \spazio \spazio 
\spazio \spazio 
\spazio  + \intTimeSpace{(\nabla \boldsymbol v)^T\boldsymbol \chi \fb{\cdot} \boldsymbol z},\\
\end{align*}
\vspace{-1cm}
\begin{align*}
{a_2^{\text{nl}}}\dual : {\Cal Y} \times {\Cal Y} \times {\Cal Y} \rightarrow {\mathbb R}  &\qquad {a_2^{\text{nl}}}\dual((\boldsymbol v, h), (\boldsymbol \chi, \lambda), (\boldsymbol \kappa, \xi)) = \displaystyle - \intSpace {h \nabla \lambda \fb{\cdot} \boldsymbol  z},\\
{a_3^{\text{nl}}}^{\ast} : {\Cal Y} \times {\Cal Y} \times {\Cal Y} \rightarrow {\mathbb R} & \qquad  {a_3^{\text{nl}}}^{\ast}((\boldsymbol v, h), (\boldsymbol \chi, \lambda), (\boldsymbol z, q)) = \displaystyle  - \intSpace {\boldsymbol v \cdot \nabla \lambda q},\\
 a_3\dual \goesto{\Cal Y}{\Cal Y}{\mathbb R} & \qquad a_3\dual ((\boldsymbol \chi, \lambda), (\boldsymbol z, q)) =   - g\intSpace {\dive( \boldsymbol \chi) q}.\\
\end{align*}
\section{Details on the Space-Time Formulation}
\label{2}
We here propose a detailed discussion of the algebraic structure briefly presented in Section \ref{FEM}. First of all, we specify the nature of the residual vector $\Cal R(\boldsymbol w \disc; \bmu)$. Then, we aim at underlining the saddle point structure of $\mathbb J(\bar w; \bmu)$. This concept is fundamental \fb{for our formulation to comply} with classical references for optimization problems such as \cite{Benzi, HinzeStokes, HinzeNS, Stoll1, Stoll}. Moreover, the saddle point structure arising from linearization justifies the reduction techniques proposed in Section \ref{aggr}, already exploited for linear steady \ocp s in \cite{negri2015reduced, negri2013reduced, rozza2012reduction} and, for time dependent problems in \cite{Strazzullo2, ZakiaMaria}.\\
To this purpose, we define $M_{\boldsymbol v}$, $M_{\boldsymbol u}$ and $M_{h}$ as mass matrices with respect to the variables $\boldsymbol v, \boldsymbol u$  and $h$, respectively, and $K$, $D$, 
$H(\boldsymbol v_k\disc)$, $\overline H\dual (\boldsymbol {v}_k \disc)$, $G(\boldsymbol v_k\disc)$, 
$G\dual (\boldsymbol {v}_k \disc)$ and $ F\dual (h_h\disc)$ 
 where: %$\overline H_{v_k} := \overline H(\boldsymbol v_k\disc)$ and  $G_{h_k} = G(h_k\disc)$,  where:
\begin{align*}
& \hspace{2mm}(K)_{ij} = a_1((\boldsymbol \phi_j, \phi_j), (\boldsymbol \phi_i, \phi_i)), \\
&\hspace{2mm}(D)_{ij} = a_2((\boldsymbol \phi_j, \phi_j), (\boldsymbol \phi_i, \phi_i))\\ 
& {(H(\boldsymbol v_k\disc))}_{ij} = a_1^{\text{nl}}((\boldsymbol v_k\disc, h_k\disc), (\boldsymbol \phi_j, \phi_j), (\boldsymbol \phi_i, \phi_i)),\\
%&  {(\overline H_{v_k})}_{ij} = a_1^{\text{nl}}( (\boldsymbol \phi_j, \phi_j), (\boldsymbol v_k\disc, h_k\disc), (\boldsymbol \phi_i, \phi_i)),
& ({ \overline H}\dual(\boldsymbol v_k\disc))_{ij} = \intSpace{(\nabla \boldsymbol v_k \disc)^T\boldsymbol \phi_i \boldsymbol \phi_j},\\
&{(G(\boldsymbol v_k\disc))}_{ij} = a_2^{\text{nl}}((\boldsymbol v_k\disc, h_k\disc), (\boldsymbol \phi_j, \phi_j), (\boldsymbol \phi_i, \phi_i)). \\
& ({G}\dual(\boldsymbol v_k\disc))_{ij} = a_3^{\text{nl}\dual}((\boldsymbol v_k \disc, h_k \disc), (\boldsymbol \phi_i, \phi_i), (\boldsymbol \phi_j, \phi_j)), \\
&({ F}\dual(h_k\disc))_{ij} = a_2^{\text{nl}\dual}((\boldsymbol v_k \disc, h_k \disc), (\boldsymbol \phi_i, \phi_i), (\boldsymbol \phi_j, \phi_j)). 
%& {(H_{h_k})}_{ij} = a_2^{\text{nl}}( (\boldsymbol \phi_j, \phi_j), (\boldsymbol v_k\disc, h_k\disc), (\boldsymbol \phi_i, \phi_i))\\
\end{align*}
Moreover, for the sake of notation, let us define the operators $S(\boldsymbol v_k\disc) = \mu_1 \Delta t K + \mu_2 \Delta t H(\boldsymbol v_k\disc)$ and  ${S}\dual(\boldsymbol v_k\disc) = \mu_1 \Delta t K - \mu_2 \Delta t H(\boldsymbol v_k\disc) 
+ \mu_2 \Delta t { \overline H}\dual(\boldsymbol v_k\disc)$. Then, the explicit form for the residual vector is
\begin{equation}
\Cal R(\boldsymbol w \disc; \bmu) = 
\underbrace{
\begin{bmatrix}
\Delta t \Cal M_{\boldsymbol v} \bar v + \Cal K_1\dual(\boldsymbol v\disc) \bar \chi + \Cal K_2\dual (h\disc) \bar \lambda \\ 
\Delta t \Cal M_{h} \bar h + \Cal K_3 \dual \bar \chi + \Cal K_4\dual(\boldsymbol v \disc) \bar \lambda \\
\alpha \Delta t \mathcal M_{\boldsymbol u}{\bar u} - \Delta t \Cal M_{\boldsymbol u} {\bar \chi} \\
\Cal K_1(\boldsymbol v\disc) \bar v + \Cal K_2 \bar h - \Delta t \Cal M_{\boldsymbol u} \bar u \\
\Cal K_4(\boldsymbol v \disc) \bar h \\
\end{bmatrix}}_{\Cal G(\bar w; \bmu)\bar w}
-
\underbrace{
\begin{bmatrix}
\Delta t \Cal M_{\boldsymbol v} \bar v_d \\
\Delta t \Cal M_{h} \bar h_d \\
\bar 0 \\
\Cal M_{\boldsymbol v}\bar v_0\\
\Cal M_{h} \bar h_0 \\
\end{bmatrix}}_{\bar f},
\end{equation}
where
$\Cal M_{\boldsymbol v}$,  $\Cal M_h$ and $\Cal M_{\boldsymbol u}$ are  block diagonal matrices with diagonal entries $\{ M_{\boldsymbol v}, \dots, M_{\boldsymbol v}\}$, $\{ M_{h}, \dots, M_{h} \}$ and $\{M_{\boldsymbol u}, \dots, M_{\boldsymbol u}\}$, respectively. 
In addition, we define the block diagonal matrices given by $\Cal K_2 = \text{diag} \{ \Delta t D, \dots, \Delta t D \} $, \newline $\Cal K_2 \dual (h\disc)= \text{diag}\{	\Delta t {  F}^{\ast}(h_1\disc), \dots, \Delta t { F}^{\ast}(h_{N_t}\disc) \}$, $\Cal K_3\dual  = \text{diag} \{	\Delta t D^{T}, \dots, \Delta t D^{T}\}$ and the matrices $\Cal K_1 (\boldsymbol v \disc)$, $\Cal K_4 (\boldsymbol v \disc)$, $\Cal K_1\dual(\boldsymbol v \disc)$ and   $\Cal K_4\dual (\boldsymbol v \disc)$, respectively, as:
\begin{equation*}\small{
\begin{bmatrix}
  M_{\boldsymbol v} +S{(\boldsymbol v_1\disc)}& 0  & \cdots & & & 0 \\ %& \Delta t D & 0 & \cdots &  &  & 0 \\ 
 % 0 &   \cdots &  &  &   \cdots& 0 \\
 -M_{\boldsymbol v} & M _{\boldsymbol v} + S{(\boldsymbol v_2\disc)}& 0  & \cdots  & & 0  \\
 0  & -M_{\boldsymbol v} & M _{\boldsymbol v} + S{(\boldsymbol v_3\disc)} &  0 & \cdots &    0 \\
& & \ddots& \ddots & &      \\ 
% &  \ddots& \ddots & &     \\ 
 0 & \cdots &  &  0 & - M_{\boldsymbol v} & M _{\boldsymbol v} + S{(\boldsymbol v_{N_t}\disc)}\\ 
 %0 &   \cdots &  &  &   \cdots& 0 \\
\end{bmatrix},}
\end{equation*}
\begin{equation*} \small{
%\begin{bmatrix}
%  \Delta t D & 0  & \cdots & & & 0 \\ %& \Delta t D & 0 & \cdots &  &  & 0 \\ 
%  M_{h} + \Delta t G_{v_1} &   \cdots &  &  &   \cdots& 0 \\
% 0 & \Delta t D & 0  & \cdots  & & 0 \\
%- M_{h}&  M_{h} + \Delta t G_{v_1}& 0 & \cdots& \cdots & 0    \\ 
%&  \ddots& \ddots & &  &   \\ 
%& &  \ddots& \ddots & &    \\ 
%0 & \cdots &  &  & 0 & \Delta t D\\ 
% 0 &   \cdots &  &  \cdots&  0& - M_{h} + \Delta t G_{v_{N_t}} \\
%\end{bmatrix},
%\end{equation*}
%\begin{equation*}
\begin{bmatrix}
  M_{h} +\Delta t G(\boldsymbol v_1 \disc) & 0  & \cdots & & & 0 \\ %& \Delta t D & 0 & \cdots &  &  & 0 \\ 
 % 0 &   \cdots &  &  &   \cdots& 0 \\
 -M_{h} & M _{h} + \Delta t G(\boldsymbol v_2 \disc)& 0  & \cdots  & & 0  \\
 0  & -M_{h} & M _{h} + \Delta t G\boldsymbol (\boldsymbol v_3 \disc) &  0 & \cdots &    0 \\
& & \ddots& \ddots & &      \\ 
% &  \ddots& \ddots & &     \\ 
 0 & \cdots &  &  0 & - M_{h} & M _{h} + \Delta t G(\boldsymbol v_{N_t} \disc)\\ 
 %0 &   \cdots &  &  &   \cdots& 0 \\
\end{bmatrix},}
\end{equation*}
\begin{equation*} \small{
\begin{bmatrix}
  M_{\boldsymbol v} +S\dual(\boldsymbol v_1\disc)& - M_{\boldsymbol v}  & \cdots & & & 0 \\ %& \Delta t D & 0 & \cdots &  &  & 0 \\ 
 % 0 &   \cdots &  &  &   \cdots& 0 \\
 0 & M _{\boldsymbol v} + S\dual(\boldsymbol v_2\disc)& - M_{\boldsymbol v}  & \cdots  & & 0  \\
 % 0  & 0 & M _{\boldsymbol v} + S_{v_3} &  - M_{\boldsymbol v} & \cdots &    0 \\
& & \ddots& \ddots & &      \\ 
% &  \ddots& \ddots & &     \\ 
0 & \cdots &  0 & M _{\boldsymbol v} + S\dual(\boldsymbol v_{{N_t}-1}\disc) & - M_{\boldsymbol v} & 0  \\  
0 & \cdots &  &   & 0 & M _{\boldsymbol v} + S\dual(\boldsymbol v_{N_t}\disc)\\ 
 %0 &   \cdots &  &  &   \cdots& 0 \\
\end{bmatrix},}
\end{equation*}
and 
\begin{equation*}\small{
\begin{bmatrix}
  M_{h} +{  G}\dual({\boldsymbol v_1 \disc}) & - M_{h}  & \cdots & & & 0 \\ %& \Delta t D & 0 & \cdots &  &  & 0 \\ 
 % 0 &   \cdots &  &  &   \cdots& 0 \\
 0 & M_{h} +{ G}\dual({\boldsymbol v_2\disc})& - M_{h}  & \cdots  & & 0  \\
 % 0  & 0 & M _{\boldsymbol v} + S_{v_3} &  - M_{\boldsymbol v} & \cdots &    0 \\
& & \ddots& \ddots & &      \\ 
% &  \ddots& \ddots & &     \\ 
0 & \cdots &  0 & M_{h} +{ G}\dual({ \boldsymbol v_{N_t - 1}\disc)} & - M_{h} & 0  \\  
0 & \cdots &  &   & 0 &  M_{h} +{ G}\dual({\boldsymbol v_{N_t}\disc)}\\ 
 %0 &   \cdots &  &  &   \cdots& 0 \\
\end{bmatrix}.}
\end{equation*}

\no Thus, the residual $\Cal R(\boldsymbol w \disc; \bmu)$ can be also written in the following compact form
\begin{equation}
\label{G}
\Cal R(\boldsymbol w \disc;  \bmu) =
\begin{bmatrix}
\Delta t \Cal M [\bar v, \bar h] + \Cal K \dual(\boldsymbol v\disc, h\disc) [\bar \chi, \bar \lambda]\\
\alpha \Delta t \Cal M_{\boldsymbol u} \bar u  - \Delta t \Cal M_{\boldsymbol u} \bar \chi\\
\Cal K(\boldsymbol v\disc) [\bar v, \bar h]  -\Delta t \Cal M_{\boldsymbol u 0} [\bar u, \bar 0] \\
\end{bmatrix}
-
\begin{bmatrix}
\Delta t \Cal M [\bar v_d, \bar h_d] \\
\bar 0 \\
\Cal M [\bar v_0, \bar h_0] \\
\end{bmatrix},
\end{equation}
where $\Cal M = 
\begin{bmatrix}
\Cal M_{\boldsymbol v}, & 0\\
0, & \Cal M_{\boldsymbol h}\\
\end{bmatrix},
$
$
\Cal K (\boldsymbol v \disc) = 
\begin{bmatrix}
\Cal K_1(\boldsymbol v \disc)& \Cal K_2 \\
0 & \Cal K_4(\boldsymbol v \disc)
\end{bmatrix}$\\
$\Cal K \dual (\boldsymbol v \disc, h \disc) = 
\begin{bmatrix}
\Cal K_1\dual(\boldsymbol v \disc)& \Cal K_2\dual(h\disc) \\
 \Cal K_3\dual & \Cal K_4 \dual(\boldsymbol v \disc)\\
\end{bmatrix}$ and 
$\Cal M_{\boldsymbol u 0} = \begin{bmatrix}
\Cal M_{\boldsymbol u}, & 0\\
0, & 0\\
\end{bmatrix}.
$ We remark that the first, the second and the last row of $\Cal R(\boldsymbol w \disc, \bmu)$ represent adjoint, optimality and state equations, respectively. We want now focus our attention on $\mathbb J(\bar w; \boldsymbol \mu)$. For the sake of clarification, we underline that with the notation $(\cdot)^\mathbb D$, we denote a quantity which \emph{derives} from the differentiation of operators. The differentiation will be applied to general space-time variables that are denoted with $\underline v, \underline h, \underline u, \underline \chi, \underline \lambda$.
Let us start our analysis with the state equation. New operators are needed: $\overline H(\boldsymbol v_k\disc)$, with ${\overline H(\boldsymbol v_k\disc)}_{ij} = a_1^{\text{nl}}((\boldsymbol \phi_j, \phi_j), (\boldsymbol v_k\disc, h_k\disc), (\boldsymbol \phi_i, \phi_i))$, and $F(h_k\disc)$ with $F(h_k\disc)_{ij} = a_2^{\text{nl}}((\boldsymbol \phi_j, \phi_j), (\boldsymbol v_k\disc, h_k\disc), (\boldsymbol \phi_i, \phi_i))$. Then, defining $$S(\boldsymbol v_k\disc)^{\mathbb D} := \fb{\mathbb D} (S(\boldsymbol v_k\disc)v_k) = \mu_1 \Delta t K + \mu_2 \Delta t H(\boldsymbol v_k\disc) + \mu_2 \Delta t \overline H(\boldsymbol v_k\disc),$$ we can differentiate the state equation as follows:
$
\mathbb D (\Cal K(\boldsymbol v \disc) \begin{bmatrix}
\bar v \\
\bar h
\end{bmatrix}  - \Delta t \Cal M_{\boldsymbol u} \bar u) $.
The process will affect only nonlinear terms and will lead to a linearized system of the form
\begin{equation}
\Cal K^{\mathbb D}(\boldsymbol v \disc, h \disc)\begin{bmatrix}
\underline v \\
\underline h
\end{bmatrix} - \Delta t \Cal M_{\boldsymbol u} \underline u,
\end{equation}
 with
$ \Cal K^{\mathbb D} =
\begin{bmatrix}
\Cal  K_1^{\mathbb D}(\boldsymbol v \disc) & \Cal K_2 \\
\Cal K_3^{\mathbb D}(h\disc) & \Cal K_4(\boldsymbol v\disc)
\end{bmatrix}
$ where $\Cal K_3^{\mathbb D} = \text{diag} \{ F(h_1\disc), \dots, F(h_{N_t}\disc) \}$ and $\Cal K_1^{\mathbb D}$ is
\begin{equation*}
\small{
\begin{bmatrix}
  M_{\boldsymbol v} + S(\boldsymbol v_1\disc)^{\mathbb D}& 0  & \cdots & & & 0 \\ %& \Delta t D & 0 & \cdots &  &  & 0 \\ 
 % 0 &   \cdots &  &  &   \cdots& 0 \\
 -M_{\boldsymbol v} & M _{\boldsymbol v} + S(\boldsymbol v_2\disc)^{\mathbb D}& 0  & \cdots  & & 0  \\
 0  & -M_{\boldsymbol v} & M _{\boldsymbol v} +  S(\boldsymbol v_3\disc) ^{\mathbb D}&  0 & \cdots &    0 \\
& & \ddots& \ddots & &      \\ 
% &  \ddots& \ddots & &     \\ 
 0 & \cdots &  &  0 & - M_{\boldsymbol v} & M _{\boldsymbol v} + S(\boldsymbol v_{N_t}\disc)^{\mathbb D}\\ 
 %0 &   \cdots &  &  &   \cdots& 0 \\
\end{bmatrix}.}
\end{equation*}
The differentiation of the optimality equation \eqref{opt_eq_k} leads to the same equation, due \fb{to} its linearity. Let us differentiate the \emph{adjoint equation}. In order to write the linearized system we define four more operators:
$H\dual(\boldsymbol {\chi}_k \disc)$, 
$\overline H \dual (\boldsymbol {\chi}_k \disc)$, $\overline F (\lambda_k \disc)$ and $\overline G(\lambda_k \disc)$ where:
\[
(H\dual(\boldsymbol {\chi}_k \disc))_{ij} =  - \intSpace{(\boldsymbol \phi_i\cdot \nabla)\boldsymbol \chi_k\disc \boldsymbol \phi_j}, \hspace{2mm}
(\overline H \dual (\boldsymbol {\chi}_k \disc))_{ij} = \intSpace{(\nabla \boldsymbol \phi_i \disc)^T\boldsymbol \chi_k\disc \boldsymbol \phi_j}, \hspace{2mm} 
\]
\[
(\overline F (\lambda_k \disc))_{ij} = a_2^{\text{nl}}((\boldsymbol \phi_j, \phi_j), (\boldsymbol \chi_k \disc, \lambda_k \disc), (\boldsymbol \phi_j, \phi_j))\]
and \[
(\overline G(\lambda_k \disc))_{ij} = a_2^{\text{nl}}((\boldsymbol \phi_j, \phi_j), (\boldsymbol \chi_k \disc, \lambda_k \disc), (\boldsymbol \phi_j, \phi_j)).
\]

\no Thanks to these quantities, we can perform the differentiation of the adjoint equation 
$$
\mathbb D(\Delta t \Cal M\begin{bmatrix}
\bar v \\
\bar h
\end{bmatrix} + \Cal K \dual(\boldsymbol v \disc,  h \disc) \begin{bmatrix}
\bar \chi \\
\bar \lambda
\end{bmatrix} ),
$$
which will result in the following linearized system
\begin{equation}
\Delta t 
\underbrace{( \Cal M + {\Cal K \dual}^{\mathbb D}(\boldsymbol \chi \disc, \lambda \disc))}_{\Cal M^{\mathbb D}(\boldsymbol \chi \disc, \lambda \disc)}
\begin{bmatrix}
\underline v \\
\underline h
\end{bmatrix} + \Cal K \dual(\boldsymbol v \disc, h \disc) \begin{bmatrix}
\underline \chi \\
\underline \lambda
\end{bmatrix},
\end{equation}
where 
$ \Cal M^{\mathbb D} = 
\begin{bmatrix}
\Cal M_1 ^{\mathbb D}(\boldsymbol \chi \disc)& \Cal M_2^{\mathbb D}(\lambda \disc) \\
\Cal M_3^{\mathbb D} (\lambda \disc)& \Cal M_4^{\mathbb D} 
\end{bmatrix}
$ with the block diagonal matrices defined by\\
$$\Cal M_1^{\mathbb D}(\boldsymbol \chi \disc)  = \text{diag}\{ M_{\boldsymbol v} + \mu_2  {H}^{\ast}(\chi_1\disc) + \mu_ 2{ \overline H}^{\ast}(\chi_1\disc)
, \dots,  M_{\boldsymbol v} + \mu_2  {H}^{\ast}(\chi_{N_t}\disc) + \mu_2 { \overline H}^{\ast}(\chi_{N_t}\disc) \},$$
$\\ \Cal M_2^{\mathbb D}(\lambda \disc) =  \text{diag}\{{ \overline F}(\lambda_{1}\disc), \dots, { \overline F}(\lambda_{N_t}\disc)\},$
$\Cal M_3^{\mathbb D}(\lambda \disc)  =  \text{diag}\{{ \overline G}(\lambda_{1}\disc), \dots, { \overline G}(\lambda_{N_t}\disc)\}$ and $\Cal M_4^{\mathbb D} = \text{diag} \{ M_{h}, \dots, M_{h}\}$.
\no We underline that, in \eqref{G} the  $\Cal K \dual \neq \Cal K^{T}$, due to the nonlinearity of the involved forms in the system, then no saddle point structure arises. However, we can recast the linearized problem in a saddle point formulation since $\Cal K \dual \equiv {\Cal K}{^{\mathbb D}}^T$. Indeed, calling with $\bar x$ the state-control space-time vector variable $
\begin{bmatrix}
\begin{bmatrix}
\bar v \\ \bar h\\
\end{bmatrix} \\
 \bar u
\end{bmatrix}$
 and with $\bar p = 
\begin{bmatrix}
\bar \chi \\ \bar \lambda
\end{bmatrix}$ the adjoint variable, and defining
\[ 
\mathsf A = 
\begin{bmatrix}
\Delta t \Cal M^{{\mathbb D}}& 0 \\
0 & \alpha \Delta t \Cal M_{\boldsymbol u} \\
\end{bmatrix}\hspace{3mm} \text{and} \hspace{3mm}
\mathsf B = 
\begin{bmatrix}
\Cal K^{\mathbb D} & -\Delta t \Cal M_{\boldsymbol u}\\
\end{bmatrix},
\]
it is simple to remark that the Frech\'et derivative can be read in the following saddle point framework:
\begin{equation}
\label{Frechet}
\mathbb J (\bar w; \boldsymbol \mu) \bar w = 
\begin{bmatrix}
\mathsf A & \mathsf B^T \\
\mathsf B & 0 \\
\end{bmatrix}
\begin{bmatrix}
\bar x \\
\bar p
\end{bmatrix}.
\end{equation}
We remark that $\mathbb J (\bar w; \boldsymbol \mu)$ is actually a generalized saddle point matrix, see \cite{Benzi} as references, where $\mathsf A \neq \mathsf A^T$. Still, we will always talk about saddle point structure from now on, since the generalization does not affect the reduced strategy used \cite{GeneralizedSaddlePoint}: indeed, the solvability condition remains the fulfillment of the inf-sup condition \cite{Babuska1971,boffi2013mixed, Brezzi74} over the state equation for the symmetric part of $\mathsf A$. Moreover, the saddle point structure does not depend on the discretization scheme used: it can be generalized for other space and time approximations.

}

\section*{Acknowledgements}
We acknowledge the support by European Union Funding for Research and Innovation -- Horizon 2020 Program -- in the framework of European Research Council Executive Agency: Consolidator Grant H2020 ERC CoG 2015 AROMA-CFD project 681447 ``Advanced Reduced Order Methods with Applications in Computational Fluid Dynamics''. We also acknowledge the INDAM-GNCS project ``Advanced intrusive and non-intrusive model order reduction techniques and applications'' and the PRIN 2017  ``Numerical Analysis for Full and Reduced Order Methods for the efficient and accurate solution of complex systems governed by Partial Differential Equations'' (NA-FROM-PDEs).
The computations in this work have been performed with RBniCS \cite{rbnics} library, developed at SISSA mathLab, which is an implementation in FEniCS \cite{fenics} of several reduced order modelling techniques; we acknowledge developers and contributors to both libraries.
\bibliographystyle{plain} % "apa":

\bibliography{BIB}

\end{document}